\documentclass[11pt,a4paper]{amsart}

\usepackage{amscd,amssymb}

\newcommand{\pfit}[1]{{\it #1:\/}}
\newenvironment{pf}{\noindent\pfit{Proof} }{\qedsymbol\par\par\medskip\par}
\makeatletter
\def\multilimits@{\bgroup
  \Let@
  \restore@math@cr
  \default@tag
 \baselineskip\fontdimen10 \scriptfont\tw@
 \advance\baselineskip\fontdimen12 \scriptfont\tw@
 \lineskip\thr@@\fontdimen8 \scriptfont\thr@@
 \lineskiplimit\lineskip
 \vbox\bgroup\ialign\bgroup\hfil$\m@th\scriptstyle{##}$\hfil\crcr}
\def\Sb{_\multilimits@}
\def\Sp{^\multilimits@}
\def\endSb{\crcr\egroup\egroup\egroup}

\makeatother

\advance\textwidth1in
\advance\oddsidemargin-0.5in
\advance\evensidemargin-0.5in
\advance\textheight1in
\advance\topmargin-0.75in

\newtheorem{theorem}[subsection]{Theorem}
\newtheorem{proposition}[subsection]{Proposition}
\newtheorem{lemma}[subsection]{Lemma}
\newtheorem{corollary}[subsection]{Corollary}

\theoremstyle{definition}

\newtheorem{definition}[subsection]{Definition}

\def\Z{\mathbb{Z}}
\def\C{\mathbb{C}}

\def\R{\mathbb{R}}
\def\Q{\mathbb{Q}}

\newcommand{\om}{\omega}
\renewcommand{\o}{\otimes}

\renewcommand{\a}{\mathcal{A}}
\renewcommand{\b}{\mathcal{B}}
\newcommand{\Braid}{\mathcal{B}\mathit{raid}}

\newcommand{\e}{\mathcal{E}}
\newcommand{\E}{\mathcal{E}}
\renewcommand{\v}{\mathcal{V}}

\renewcommand{\r}{\mathcal{R}}
\newcommand{\w}{\mathcal{W}}

\renewcommand{\SS}{\mathbb{S}}
\newcommand{\bull}{{\bullet}}

\newcommand{\CO}{\mathcal{O}}

\newcommand{\1}{{1\!\!1}}

\newcommand{\Hom}{\operatorname{Hom}}

\newcommand{\abuts}{\Rightarrow}

\newcommand{\Aut}{\operatorname{Aut}}
\newcommand{\Aff}{\operatorname{Aff}}
\newcommand{\PSL}{\operatorname{PSL}}

\newcommand{\Ind}{\operatorname{Ind}}
\newcommand{\Res}{\operatorname{Res}}
\renewcommand{\span}{\operatorname{span}}
\newcommand{\p}{\partial}

\newcommand{\CM}{\mathcal{M}}

\newcommand{\Schur}{\mathsf{S}}

\newcommand{\CL}{\mathcal{L}}

\newcommand{\gr}{\operatorname{gr}}

\renewcommand{\*}{\cdot}
\renewcommand{\L}{\mathsf{L}}
\newcommand{\TT}{\mathbb{T}}

\def\<{\langle}
\def\>{\rangle}
\def\Flag{\operatorname{Flag}}

\def\Ass{{\operatorname{\mathcal{A}\mathit{ss}}}}
\def\Com{{\operatorname{\mathcal{C}\mathit{om}}}}
\def\Lie{{\operatorname{\mathcal{L}\mathit{ie}}}}
\def\sgn{\operatorname{sgn}}

\newcommand{\Mbar}{\overline{\mathcal{M}}}
\newcommand{\Grav}{\mathcal{G}\mathit{rav}}
\def\mbar{\mathcal{H}\mathit{ycom}}
\def\m{\mathbf{m}}

\def\.{\wedge}
\def\({(\!(}
\def\){)\!)}

\def\VV{\operatorname{Vert}}
\def\Edge{\operatorname{Edge}}
\def\Leg{\operatorname{Leg}}

\def\]{{]\!]}}
\def\[{{[\![}}
\def\ch{\operatorname{ch}}
\def\Ch{\operatorname{Ch}}
\def\BB{\mathsf{B}}
\def\BC{\mathcal{B}}
\def\CP{\mathbb{CP}}
\newcommand{\Dual}{\vee}

\newcommand{\Tr}{\operatorname{Tr}}
\newcommand{\Exp}{\operatorname{Exp}}
\newcommand{\Log}{\operatorname{Log}}
\newcommand{\CT}{\mathcal{T}}
\newcommand{\Wedge}{\Lambda}

\begin{document}

\title{Operads and moduli spaces of genus $0$ Riemann surfaces}

\author{E. Getzler}

\address{Department of Mathematics, MIT, Cambridge MA 02139 USA}

\email{getzler@@math.mit.edu}

\thanks{The author is partially supported by a fellowship of the Sloan
Foundation and a research grant of the NSF.}

\maketitle

In this paper, we study two dg (differential graded) operads related to the
homology of moduli spaces of pointed algebraic curves of genus $0$. These
two operads are dual to each other, in the sense of Kontsevich
\cite{Kontsevich} and Ginzburg and Kapranov \cite{GK}.

This duality is analogous to the duality between commutative and Lie
algebras, which goes back to Quillen \cite{Quillen}. Let us recall what form
this duality takes. Note that all commutative algebras in this paper will
be non-unital --- that is, we do not take the existence of an identity as
one of the axioms.

Associate to a dg commutative algebra $A$ the chain complex $\L QA$ whose
underlying graded vector space is the desuspension of the free graded Lie
coalgebra generated by $\Sigma A$, and whose differential is the sum of the
internal differential of $A$ and the coderivation induced by the product of
$A$, thought of as a map from $\Sigma^{-1}(\Sigma A\.\Sigma A)$ to
$\Sigma^{-1}(\Sigma A)$. (Thus, $\L QA$ is the desuspenion of the Harrison
complex of $A$.) There is a natural map from $\L QA$ to the complex of
irreducibles $QA=A/A^2$ of $A$.

Similarly, associate to a dg Lie algebra $L$ the chain complex $\L QL$
whose underlying graded vector space is the desuspension of the free graded
commutative coalgebra generated by $\Sigma L$, and whose differential is
the sum of the internal differential of $L$ and the coderivation induced by
the bracket of $L$, thought of as a map from $\Sigma^{-1}(\Sigma
L\vee\Sigma L)$ to $\Sigma^{-1}(\Sigma L)$. (Thus, $\L QA$ is the
desuspenion of the Chevalley-Eilenberg complex of $L$.) There is a natural
map from $\L QL$ to the complex of irreducibles $QL=L/[L,L]$ of $L$.

Over a field of characteristic zero, the following basic result holds; the
first part is due to Koszul, and the second part is due to Barr and
Quillen. (In fact, the first part holds in any characteristic, but we will
restrict attention to fields of characteristic zero throughout this paper.)
\begin{theorem} \label{Harrison}
If $A$ is a free dg commutative algebra, the map $\begin{CD}\L QA@>>>QA
\end{CD}$ is a
homotopy equivalence (that is, a quasi-isomorphism of chain complexes). If
$L$ is a free dg Lie algebra, the map $\begin{CD}\L QL@>>>QL
\end{CD}$ is a homotopy
equivalence.
\end{theorem}

In this paper, we prove an analogue of this theorem in which the r\^ole of
commutative algebras is taken by a beautiful algebraic structure discovered
by Dijkgraaf, Verlinde and Verlinde \cite{DVV}, which we call a
hypercommutative algebra. This is a chain complex $A$ with a sequence of
graded symmetric products $\begin{CD}(x_1,\dots,x_n):A^{\o n}@>>>A
\end{CD}$ of degree
$2(n-2)$, which satisfy the following generalized associativity condition:
if $a,b,c,x_1,\dots,x_n$, $n\ge0$, are elements of $A$,
\begin{equation} \label{hypercommutative}
\sum_{S_1\amalg S_2=\{1,\dots,n\}} \pm ((a,b,x_{S_1}),c,x_{S_2})
= \sum_{S_1\amalg S_2=\{1,\dots,n\}} \pm (a,(b,c,x_{S_1}),x_{S_2}) .
\end{equation}
Here, if $S=\{s_1,\dots,s_k\}$ is a finite set, $x_S$ is an abbreviation
for $x_{s_1},\dots,x_{s_k}$. The symbol $\pm$ indicates the
Quillen sign convention for $\Z/2$-graded vector spaces: it equals $+1$ if
all of the variables have even degree.

In the cases $n=0$ and $1$, we obtain respectively the relations
$(a,(b,c))=((a,b),c)$ and
$$
(a,(b,c,d)) + (a,(b,c),d)
= ((a,b),c,d) + (-1)^{|c|\,|d|} ((a,b,d),c) .
$$
In particular, the product $(a,b)$ is a graded commutative, associative
product in the usual sense.

We may think of the products in a hypercommutative algebra as the Taylor
coefficients of a formal deformation of the commutative product $(a,b)$,
parametrized by the double suspension $\Sigma^{-2}A$ of $A$:
$$
(a,b)_t = \sum_{k=0}^\infty \frac{1}{k!}
(a,b,\underset{\text{$k$ times}}{\underbrace{t,\dots,t}}) .
$$
The associativity of $(a,b)_t$ is equivalent to the sequence of identities
satisfied by the products $(x_1,\dots,x_n)$.

Recently, Kontsevich and Manin have shown that under certain circumstances,
the quantum (or Floer) cohomology of a compact K\"ahler manifold is a
hypercommutative algebra \cite{KM}; this shows that the rather heuristic
calculations of Dijkgraaf, Verlinde and Verlinde may be justified by the
study of the moduli spaces $\Mbar_{0,n}$. We refer to this reference for
further details; extensions to a larger class of symplectic manifolds may
be found in \cite{MS} and \cite{RT}.

In our analogue of Theorem \ref{Harrison}, commutative algebras are
replaced by hypercommutative algebras, while Lie algebras are replaced by
gravity algebras, introduced by the author in \cite{weil}. A gravity
algebra is a chain complex with graded antisymmetric products
$\begin{CD}[x_1,\dots,x_n]:A^{\o n}@>>>A
\end{CD}$ of degree $2-n$, satisfying the following
relations: if $k>2$ and $\ell\ge0$, and $a_1,\dots,a_k,b_1,\dots,b_\ell\in
A$,
\begin{equation} \label{gravity}
\sum_{1\le i<j\le k}
\pm [[a_i,a_j],a_1,\dots,\widehat{a_i},\dots,\widehat{a_j},\dots,a_k,
b_1,\dots,b_\ell]
= \begin{cases} [[a_1,\dots,a_k],b_1,\dots,b_\ell] , & \ell>0 , \\
0 , & \ell=0 . \end{cases}
\end{equation}
For example, setting $k=3$ and $\ell=0$, we obtain the Jacobi relation for
$[a,b]$.

Our proof of the duality between hypercommutative and gravity algebras is
based on work of Beilinson and Ginzburg \cite{BG} and Ginzburg and Kapranov
\cite{GK}. If $X$ is a topological space, let $X^n_0$ be the configuration
space of $n$ distinct labelled points
$$
X^n_0 = \text{embeddings of $\{1,\dots,n\}$ in $X$} .
$$
Let $\CM_{0,n}=(\CP^1)^n_0/\PSL(2,\C)$, $n\ge3$, be the moduli space of
smooth projective curve $\Sigma$ of genus $0$ with $n$ marked points: it is
a smooth variety of dimension $n-3$.

The moduli spaces $\CM_{0,n}$ have smooth compactifications $\Mbar_{0,n}$,
which were constructed by Deligne, Mumford and Knudsen \cite{Knudsen}. The
manifold $\Mbar_{0,n}$ is a stratified space, whose closed strata
$\Mbar\(T\)$ are labelled by trees $T$ with $n$ legs. The above authors
studied the duality between commutative and Lie algebras by applying mixed
Hodge theory to the spectral sequence for the inclusion
$\begin{CD}\CM_{0,n}@>>>\Mbar_{0,n}
\end{CD}$. Keel \cite{Keel} has shown that the cycles
$[\Mbar\(T\)]$ carried by these strata span the homology of $\Mbar_{0,n}$,
and building on his work, Kontsevich and Manin \cite{KM} have found a
complete set of relations among these cycles; in fact, these relations are
in a sense exactly the relations which hold in a hypercommutative
algebra. We will recover their results using mixed Hodge theory.

This paper is written in the language of operads, and may be considered to
be an appendix to \cite{GK}. We extend their notion of Koszul operads to
operads which are generated by operations with more than two operands, such
as the hypercommutative and gravity operads. (The idea of generalizing
Koszul duality to this setting came from conversations with Kontsevich. I
wish to thank him for their permission to present those ideas here.)

As a byproduct of our work, we obtain new formulas for the characters of
the $\SS_n$-modules $H_i(\CM_{0,n})$ and $H_i(\Mbar_{0,n})$. These formulas
illustrate the use of the Legendre transform for symmetric functions
introduced in \cite{modular}. Of course, it would be very interesting if
any of the ideas of the paper applied to higher genus moduli spaces of
curves. However, since the purity of the mixed Hodge structure of
$\CM_{0,n}$ plays such a central r\^ole here, and is in some sense equivalent
to our main theorem, any extension to higher genus will be rather subtle.

The research of the author is partially supported by the NSF and the
A.P. Sloan Research Foundation. The author is grateful to J.D.S. Jones,
M.M. Kapranov and particularly M. Kontsevich for the many ways in which they
assisted in the preparation of this article.

\section{$\SS$-modules and operads}

In this section, we recall the basic definitions of the theory of
operads. For more details, see \cite{n-algebras} and \cite{GK}.

\subsection{$\SS$-modules}
An $\SS$-module is a collection of chain complexes (all chain complexes in
this paper are over the field $\C$, and have finite dimensional total
homology)
$$
\{ \v(n) \mid n\ge0 \} ,
$$
together with an action of $\SS_n$ on $\v(n)$. This definition generalizes
Joyal's notion of a linear species \cite{Joyal}, which is an ungraded
$\SS$-module.

A chain complex $V$ may be thought of as an $\SS$-module by setting
$$
\v(n) = \begin{cases} V , & n=1 , \\ 0 , & n\ne1 . \end{cases}
$$

\subsection{Schur functors}
Given an $\SS$-module $\v$ and a finite set $S$, we define
$$
\v(S) = \Bigl(
\bigoplus\begin{Sb}\text{bijections}\\\begin{CD}f:\{1,\dots,n\}@>>>S
\end{CD}\end{Sb}
\v(n) \Bigr)_{\SS_n} .
$$
It is clear that if $S=\{1,\dots,n\}$, then $\v(S)$ is naturally
identified with $\v(n)$.

To an $\SS$-module $\v$ is associated an endofunctor of the category of
chain complexes, called the Schur functor of $\v$, by the formula
$$
V \mapsto \Schur(\v,V) = \bigoplus_{n=0}^\infty \v\o_{\SS_n} V^{\o n} ;
$$
here $V^{\o n}$ is the graded $n$th tensor power of $V$. Introduce a
monoidal structure on the category of $\SS$-modules, with tensor product
$$
(\v\circ\w)(n) = \bigoplus_{k=0}^\infty \Bigl( \v(k) \o
\bigoplus_{\begin{CD}f:\{1,\dots,n\}@>>>\{1,\dots,k\}
\end{CD}} \bigotimes_{i=1}^k
\w(f^{-1}(i)) \Bigr)_{\SS_k} ,
$$
and unit the $\SS$-module $\1$:
$$
\1(n) = \begin{cases} \C , & n=1 , \\ 0 , & n\ne1 .
\end{cases}$$
The peculiar formula for $\v\circ\w$ is justified by
$$
\Schur(\v,\Schur(\w,V)) = \Schur(\v\circ\w,V) .
$$
Note also that $\Schur(\1,V)=V$. For more on this formalism, see Chapter 1
of \cite{n-algebras}.

\subsection{Operads}
An operad is a monoid in the category of $\SS$-modules, that is, an
$\SS$-module $\a$ with product $\begin{CD}\rho:\a\circ\a@>>>\a
\end{CD}$ and unit
$\begin{CD}\eta:\1@>>>\a
\end{CD}$ satisfying the axioms of associativity and unit
\cite{Maclane}. We denote the image under the product $\rho$ of
$$
a\o b_1\o\dots\o b_k \in \a(k)\o\a(n_1)\o\dots\o\a(n_k)
$$
by $a(b_1,\dots,b_k)$ and the unit by $1\in\a(1)$.

An operad structure on an $\SS$-module $\a$ such that $\a(n)=0$ for $n\ne1$
is the same thing as an associative algebra structure on $\a(1)$. Whereas
an element of an algebra has only one ``input'' and one ``output,'' an
element of an operad has multiple inputs and one output.

\subsection{The endomorphism operad of a chain complex}
If $V$ is a chain complex, its endomophism operad is the $\SS$-module
$$
\e_V(n) = \Hom(V^{\o n},V) .
$$
This is an operad, whose product is given by composition: if
$a\in\Hom(V^{\o k},V)$ and $b_i\in\Hom(V^{\o n_i},V)$, then
$$
a(b_1,\dots,b_k) = a \* (b_1\o\dots\o b_k) ,
$$
where we think of $b_1\o\dots\o b_k$ as an element of
$\Hom(V^{\o(n_1+\dots+n_k)},V^{\o k})$.

\subsection{Suspension of operads}
If $V$ is a chain complex, denote by $\Sigma V$ the chain complex such that
$(\Sigma V)_i=V_{i-1}$, with differential $-\delta$. If $\v$ is an
$\SS$-module, denote by $\Lambda\v$ the $\SS$-module
$$
(\Lambda\v)(n) = \Sigma^{1-n} \sgn_n \o \v(n) ,
$$
where $\sgn_n$ is the sign character of $\SS_n$. There is a natural
isomorphism
$$
\Schur(\Lambda\v,V) \cong \Sigma\Schur(\v,\Sigma^{-1}V) .
$$
It follows that if $\a$ is an operad, then so is $\Lambda\a$, and if $A$
is an $\a$-algebra, then $\Sigma A$ is a $\Lambda\a$-algebra.

\subsection{Algebras over an operad}
An algebra over an operad $\a$ is a chain complex $A$, together with a
morphism of operads $\begin{CD}\a@>>>\e_A
\end{CD}$. Thus, if $A$ is an algebra over an operad
$\a$ and $\rho\in\a(n)$, there is a product $a_1\o\dots
a_n\mapsto\rho(a_1,\dots,a_n)$ from $A^{\o n}$ to $A$. These products are
equivariant, under the actions of $\SS_n$ on $\a(n)$ and $A^{\o n}$,
associative with respect to the product of $\a$, and $1(a)=a$, where
$1\in\a(1)$ is the unit of $\a$.

\subsection{Construction of operads}
Given an algebraic structure defined by a set of multilinear operations
together with a set of multilinear relations among them, one may construct
the operad $\a$ having this presentation, in such a way that an $\a$-algebra
is the same thing as an instance of the original algebraic structure.

Denote by $V\mapsto T(V)$ the free algebra generated by the chain complex
$V$ with respect to this algebraic structure. To define $\a(n)$, we form
the free algebra $T(x_1,\dots,x_n)$ generated by the free vector space
$\C^n$ of rank $n$. The torus $(\C^\times)^n$ acts on this chain complex;
let $\a(n)$ be the $\SS_n$-submodule on which it acts by the character
$(z_1,\dots,z_n)\mapsto z_1\dots z_n$. The group $\SS_n$ acts on
$T(x_1,\dots,x_n)$, and thus on $\a(n)$, by permutation of the letters
$x_i$. It is not difficult to see that the $\SS$-module thus constructed is
an operad: the unit is the word $x_1\in\a(1)\subset T(x_1)$, while the
product is defined by substitution.

Let us give some examples of this construction. In the case of commutative
algebras, we call the resulting operad $\Com^+$; the $\SS_n$-module
$\Com^+(n)$ is spanned by the word $x_1\dots x_n$ the free commutative
algebra generated by letters $\{x_1,\dots,x_n\}$, and carries the trivial
action of $\SS_n$. (Here, we are dealing with non-unital commutative
algebras, so that $\Com^+(0)=0$.)

Let $\Ass^+$ be the operad associated to associative algebras; the
$\SS_n$-module $\Ass^+(n)$ is spanned by the words
$$
\{ x_{\sigma(1)}\dots x_{\sigma(n)} \mid \sigma\in\SS_n \} ,
$$
and carries the regular representation of $\SS_n$. (Again, we set
$\Ass^+(0)=0$.)

Finally, the operad associated to Lie algebra is denoted $\Lie^+$. The
underlying $\SS$-module may be studied by means of the
Poincar\'e-Birkhoff-Witt theorem, which implies that
$$
\Ass^+ \cong \Com^+ \circ \Lie^+ .
$$
In Proposition \ref{Lie}, we will show how this leads to a formula for the
character of $\Lie^+(n)$, due to Klyachko \cite{Klyachko}.

\subsection{The configuration spaces $\C^n_0$ and the braid operad $\b$}
Let $\C^n_0$ be the configuration space of $n$ labelled points in
$\C$. Define the $\SS$-module $\Braid$ by
$$
\Braid(n) = \begin{cases} H_\bull(\C^n_0) , & n>0 , \\ 0 , & n=0 .
\end{cases}$$
We now construct a natural operad structure on $\Braid$. (We called this the
braid operad in \cite{n-algebras}.)

In the definition of operads, one can replace the category of chain
complexes by the category of topological spaces, and the tensor product by
Cartesian product, obtaining the notion of a topological operad. If $\CO$
is a topological operad and $H_\bull(-)$ is a generalized homology theory
with products, $H_\bull(\CO)$ is an operad in the category of graded vector
spaces; this gives a useful method of constructing operads.

Boardman and Vogt have constructed a topological operad called the little
discs operad \cite{BoardmanVogt}. Let $D$ be the closed unit disc in $\C$,
and let $\CO(n)$ be the topological space
$$\textstyle
\CO(n) = \bigl\{ \binom{z_1\dots z_n}{r_1\dots r_n}
\in \binom{D^n}{\R^n_+} \big|
\text{ the discs $r_iD+z_i$ are disjoint subsets of $D$} \bigr\} .
$$
The symmetric group $\SS_n$ acts on $\CP(n)$ by permuting the discs:
$$\textstyle
\sigma\binom{z_1\dots z_n}{r_1\dots r_n} =
\binom{z_{\sigma(1)}\dots z_{\sigma(n)}}{r_{\sigma(1)}\dots r_{\sigma(n)}} .
$$
The product in this operad is defined by gluing of disks: if
$a=\binom{z_1\dots z_k}{r_1\dots r_k}$ and
$b_i=\binom{y_{i,1}\dots y_{i,n_i}}{s_{i,1}\dots s_{i,n_i}}$, then
$$\textstyle
a(b_1,\dots,b_k) = \begin{pmatrix}
r_1y_{1,1}+z_1&\dots&r_1y_{1,n_1}+z_1&\dots&
r_ky_{k,1}+z_k&\dots&r_ky_{k,n_k}+z_k \\
r_1s_{1,1}&\dots&r_1s_{1,n_1}&\dots&r_ks_{k,1}&\dots&r_ks_{k,n_k}
\end{pmatrix} .
$$
The map $\begin{CD}\CO(n)@>>>\C^n_0
\end{CD}$ defined by
$\binom{z_1\dots z_n}{r_1\dots r_n}\mapsto(z_1,\dots,z_n)$ is a homotopy
equivalence, and thus the homology operad $H_\bull(\CO)$ of this
topological operad has $\Braid$ as its underlying $\SS$-module.

The operad $\Braid$ has the following presentation (see \cite{Cohen:thesis}
and \cite{n-algebras}): it is generated by two operations, a commutative
product of degree $0$ and a Lie bracket of degree $1$, satisfying the
Poisson relation:
$$
[a,bc] = [a,b]c + (-1)^{(|a|+1)|b|} b[a,c] .
$$
It follows that the $\SS$-module $\Braid$ is isomorphic to
$\Com^+\circ\Lambda^{-1}\Lie^+$. In particular, we obtain another
realization of the $\SS_n$-module $\Lie^+(n)$:
$$
\Lie^+(n) \cong \sgn_n \o H_{n-1}(\C^n_0) .
$$

\section{Cyclic operads}

All of the operads which we discuss in this article, notably the
hypercommutative and gravity operads, are cyclic operads
\cite{cyclic} --- that is, there is a notion of invariant inner product on
algebras over these operads. In this section, we recall the basics of the
theory of cyclic operads.

\subsection{Cyclic $\SS$-modules}
A cyclic $\SS$-module is an $\SS$-module $\v$ together with an action of
$\SS_{n+1}$ on $\v(n)$ extending the action of $\SS_n$. The name derives
from the fact that an action of $\SS_{n+1}$ is determined by compatible
actions of $\SS_n$ and the cyclic group $C_{n+1}\subset\SS_{n+1}$ generated
by the cycle $(01\dots n)$. We denote the action of $(01\dots n)$ on
$\v(n)$ by $v\mapsto v^*$, motivated by the fact that a cyclic $\SS$-module
structure on the $\SS$-module associated to a chain complex $V$ is just an
involution $v\mapsto v^*$. It is convenient to write $\v\(n\)$ for
$\v(n+1)$.

If $\v$ is a cyclic $\SS$-module, denote by $\Lambda\v$ the cyclic
$\SS$-module
$$
(\Lambda\v)\(n\) = \Sigma^{2-n} \sgn_n \o \v\(n\) .
$$
Thus, $\Lambda$ applied to the $\SS$-module underlying $\v$ is isomorphic
to the $\SS$-module underlying $\Lambda\v$.

A stable cyclic $\SS$-modules is a cyclic $\SS$-module $\a$ which satisfies
the condition that $\a\(n\)=0$ for $n<3$. (The word ``stable'' comes from
the theory of algebraic curves.)

\subsection{Cyclic operads}
A cyclic operad $\a$ is a cyclic $\SS$-module $\a$ with an operad
structure, such that $1^*=1$ and for all $a\in\a(k)$ and $b\in\a(l)$,
$$
a(1,\dots,1,b)^* = (-1)^{|a|\,|b|} b^*(a^*,1,\dots,1) .
$$
Note that if $\v$ is a cyclic $\SS$-module associated to a chain complex
$V$ with involution, then a cyclic operad structure on $\v$ is the same
thing as a $\ast$-algebra structure on $V$.

\subsection{Stable cyclic operads}
A cyclic operad $\a$ whose underlying $\SS$-module is stable is called a
stable cyclic operad. Because of this condition, $\a\(2\)=0$, so we can no
longer think of $\a$ as having a unit. This requires the introduction of
non-unital operads, following Markl. Observe that all of the products
$a(b_1,\dots,b_k)$ of a (unital) operad may be obtained by iterating the
products
$$
a\circ_i b = a(\underset{\text{$i-1$ times}}{\underbrace{1,\dots,1}},b,
\underset{\text{$k-i$ times}}{\underbrace{1,\dots,1}}) .
$$
The axioms for non-unital operads may be found in \cite{modular}: they are
of two types, equivariance and associativity. Finally, a non-unital operad
is cyclic if $(a\circ_kb)^*=(-1)^{|a|\,|b|}b^*\circ_1a^*$.

\subsection{Invariant inner products and cyclic algebras}
If $V$ is a chain complex with inner product $\<-,-\>$ (which we suppose to
be non-degenerate), the endomorphism operad $\e_V$ is a cyclic operad,
such that if $a\in\a(n)$ and $v_i\in V$, $0\le i\le n$,
$$
\< v_0 , a(v_1,\dots,v_n) \> = \< v_n , a^*(v_0,\dots,v_{n-1}) \> .
$$
A cyclic algebra over a cyclic operad is a chain complex $A$ with inner
product, together with a morphism of cyclic operads $\begin{CD}\a@>>>\e_A
\end{CD}$.

If $A$ is an algebra over the operad underlying a cyclic operad $\a$, and
$\<-,-\>$ is an inner product on $A$, we say that the inner product is
invariant.

Note that if $\a$ is a cyclic operad, then the operad $\Lambda\a$ is not
cyclic, but rather anticyclic:
$(a\circ_kb)^*=-(-1)^{|a|\,|b|}b^*\circ_1a^*$. This reflects the fact that
if $V$ is a chain complex and $\<-,-\>$ is an inner product on $V$, then
there is induced on $\Sigma V$ an antisymmetric non-degenerate bilinear
form $(-1)^{|v|}\<\Sigma v,\Sigma w\>$. However, $\Lambda^2\a$ is again a
cyclic operad, and if $A$ is a cyclic $\a$-algebra, then $\Sigma^2A$ is a
cyclic $\Lambda^2\a$-algebra.

\subsection{Examples of cyclic operads}
Since associative, commutative and Lie algebras all have well-known notions
of invariant inner product, it is not surprising that the corresponding
operads are cyclic. For the operads $\a^+$, where
$\a\in\{\Ass,\Com,\Lie\}$, define the $\SS$-module $\a$ by setting
$\a(n)=\a^+(n)$ for $n>1$ and $\a(n)=0$ for $n\le1$.

The cyclic structure of $\Com$ is simple to describe, since an inner
product on a commutative algebra is invariant if and only if
\begin{equation} \label{invariant}
\<a,bc\> = \<ab,c\> .
\end{equation}
It follows from this formula that $a=a^*$ for all $a\in\Com\(n\)$.

In a similar way, the action of $\SS_n$ on $\Ass\(n\)$ is determined by the
condition that an inner product on an associative algebra is invariant if
and only if \eqref{invariant} holds. It turns out that the $\SS_n$-module
$\Ass\(n\)$ is the induced representation $\Ind_{C_n}^{\SS_n}\1$, where
$C_n$ is the subgroup of $\SS_n$ generated by $\tau_n$, and $\1$ is its
trivial representation. We may think of $\Ass\(n\)$ as being the
$\SS_n$-module spanned by symbols
$$
\<x_{\sigma(1)},x_{\sigma(2)}\dots x_{\sigma(n)}\> , \quad
\sigma\in\SS_n ,
$$
representing the inner product of $x_{\sigma(2)}\dots x_{\sigma(n)}$ with
$x_{\sigma(1)}$, subject to the relations
$$
\<x_{\sigma(1)},x_{\sigma(2)}\dots x_{\sigma(n)}\> \sim
\<x_{\sigma(n)},x_{\sigma(1)}\dots x_{\sigma(n-1)}\> .
$$

The cyclic structure on the Lie operad is associated with the usual notion
of an invariant inner product (or Killing form) on a Lie algebra,
satisfying
$$
\<[a,b],c\> = \<a,[b,c]\> .
$$
The character of the $\SS_{n+1}$-module $\Lie(n)$ is calculated in
\cite{modular}. We will obtain a realization of this representation as a
homology group in the next section.

It is proved in \cite{cyclic} that the braid operad $\Braid$ is not a
cyclic operad.

\section{The moduli spaces $\Mbar_{0,n}$}

In this section, we study the combinatorial structure of the compactified
moduli spaces $\Mbar_{0,n}$. We then define the gravity and
hypercommutative operads, and introduce the fundamental exact sequences
relating them, which are obtained by considering the spectral sequence
associated by Deligne to the stratified space $\Mbar_{0,n}$.

\subsection{Graphs and trees}
The strata of the compactification $\Mbar_{0,n}$ are labelled by trees with
$n$ legs, and we recall some definitions from the theory of trees in this
paragraph (see also \cite{GK} and \cite{modular}).

A graph $G=(F,\pi,\tau)$ is a finite set $F=\Flag(G)$, the set of flags of
the graph, together with a partition $\pi$ and an involution $\tau$ of
$F$. (By a partition, we mean a decomposition of $F$ into disjoint subsets,
possibly empty, called its blocks.)

The vertices $\VV(G)$ of the graph $G$ are the blocks of $\pi$, the edges
$\Edge(G)$ are the $2$-cycles of $\tau$, while the legs $\Leg(G)$ are the
fixed points of $\tau$. To a graph $G$ is associated a cell complex $G$
whose cells have dimension $0$ and $1$, and whose ends correspond to the
legs $\Leg(G)$. A graph is called a tree if this complex is simply
connected. We will have no further use for non-simply connected graphs in
this paper; however, much of the theory we describe has an analogue for
general graphs \cite{modular}.

The legs $\Leg(v)$ of a vertex $v\in\VV(G)$ are the flags in the
corresponding equivalence class, while the valence $|v|$ of a vertex is the
cardinality of $\Leg(v)$.

If $S$ is a finite set, let $\CT\(S\)$ be the set of isomorphism classes of
trees $T$ whose external edges are labelled by the elements of $S$ and such
that each vertex has valence at least three. Note that $\CT\(S\)$ is
finite. The set of trees is graded by the number of edges:
$$
\CT\(S\) = \bigcup_{i=0}^{|S|-3} \CT_i\(S\) .
$$
In particular, $\CT_0\(S\)$ has a single element, the tree with one vertex
whose set of flags equals $S$.

Denote by $\det(S)$ the determinant line $\Wedge^{\text{max}}\C^S$, which
is a representation of $\Aut(S)$. (For example, $\det(\{1,\dots,n\})$ is
just the sign representation $\sgn_n$.) If $T$ is a tree, let $\det(T)$ be
the determinant line $\det(\VV(T))$ of the set of vertices of $T$. There
are natural isomorphisms
\begin{equation} \label{det}
\det(T) \cong \det(\Edge(T))
\cong \det(\Leg(T)) \o \bigotimes_{v\in\VV(T)} \det(\Leg(v)) .
\end{equation}

\subsection{Stable curves}
A stable curve with $n$ marked points is a projective curve $\Sigma$ whose
only singularities are double points, together with an embedding of
$\{1,\dots,n\}$ in the set of smooth points of $\Sigma$, such that there
are no continuous automorphisms of $\Sigma$ fixing the marked points and
double points. Knudsen has proved that the moduli space $\Mbar_{g,n}$ of
stable curves of arithmetic genus $g$ with $n$ marked points is a compact
orbifold, obtained by adjoining to $\CM_{g,n}$ a divisor with normal
crossings \cite{Knudsen}.  We will only be interested in the genus zero
case, in which case $\Mbar_{0,n}$ is actually a projective variety. Another
reference for the genus $0$ case is \cite{Kapranov}.

A stable curve $\Sigma$ determines a graph $\Gamma_\Sigma$, called the dual
graph of $\Sigma$: the vertices of $\Gamma_\Sigma$ are the components of
the subvariety of smooth points of $\Sigma$, the edges are the double
points, and $\Leg(\Gamma_\Sigma)$ is the set of marked points. In
particular, if $\Sigma$ has arithmetic genus $0$ and $n$ marked points,
then $\Gamma_\Sigma\in\CT\(n\)$.

The moduli space $\Mbar_{0,n}$ is a stratified space: it has one stratum
$\CM\(T\)$ for each tree $T\in\CT\(n\)$, consisting of all curves
$\Sigma\in\Mbar_{0,n}$ such that $\Gamma_\Sigma=T$. The stratum $\CM\(T\)$
is isomorphic to the product
$$
\prod_{v\in\VV(T)} \CM_{0,|v|} ,
$$
and has codimension equal to the number of edges of $T$. For example
$\Mbar_{0,4}$ has the following four strata:
$$
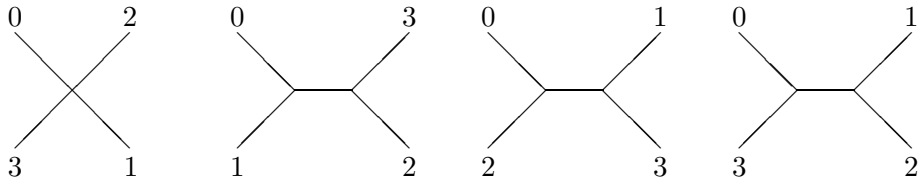
\setlength{\unitlength}{0.0075in}%
\begin{picture}(800,130)(25,0)
\put(100,5){\begin{picture}(80,102)(80,650)
\put( 80,740){\line( 1,-1){ 80}}
\put(160,740){\line(-1,-1){ 80}}
\put( 75,745){0}
\put(155,745){2}
\put( 75,640){3}
\put(156,640){1}
\end{picture}}
\put(250,0){\begin{picture}(125,107)(55,605)
\put( 60,700){\line( 1,-1){ 40}}
\put(100,660){\line(-1,-1){ 40}}
\put(100,660){\line( 1, 0){ 40}}
\put(140,660){\line( 1, 1){ 40}}
\put(140,660){\line( 1,-1){ 40}}
\put( 55,705){0}
\put(175,705){3}
\put( 55,600){1}
\put(175,600){2}
\end{picture}}
\put(425,0){\begin{picture}(125,107)(55,605)
\put( 60,700){\line( 1,-1){ 40}}
\put(100,660){\line(-1,-1){ 40}}
\put(100,660){\line( 1, 0){ 40}}
\put(140,660){\line( 1, 1){ 40}}
\put(140,660){\line( 1,-1){ 40}}
\put( 55,705){0}
\put(175,705){1}
\put( 55,600){2}
\put(175,600){3}
\end{picture}}
\put(600,0){\begin{picture}(125,107)(55,605)
\put( 60,700){\line( 1,-1){ 40}}
\put(100,660){\line(-1,-1){ 40}}
\put(100,660){\line( 1, 0){ 40}}
\put(140,660){\line( 1, 1){ 40}}
\put(140,660){\line( 1,-1){ 40}}
\put( 55,705){0}
\put(175,705){1}
\put( 55,600){3}
\put(175,600){2}
\end{picture}}
\end{picture}$$

If $T\in\CT_k(n)$, denote by $\Mbar\(T\)$ the closure of the stratum
$\CM\(T\)$ of $\Mbar_{0,n}$, and by $[\Mbar\(T\)]$ the corresponding cycle
in $H_{2(n-k-3)}(\Mbar_{0,n})$. The following theorem is due to Keel
\cite{Keel}.
\begin{theorem} \label{Keel:1}
The cycles $[\Mbar\(T\)]$, $T\in\CT\(n\)$ span $H_\bull(\Mbar_{0,n})$.
\end{theorem}

In Section \ref{Keel:proof}, we give a proof of this theorem which differs
from Keel's, and uses mixed Hodge theory.

\subsection{The gravity operad}
Let $\Grav$ be the stable cyclic $\SS$-module
$$
\Grav\(n\) = \begin{cases}
\Sigma^{3-n}\sgn_n\o H_\bull(\CM_{0,n}) , & n\ge3 , \\ 0 , & n<3 .
\end{cases}$$
Note that $\Grav\(n\)$ is concentrated in degrees $3-n\le i\le 0$.

There is a natural cyclic operad structure on $\Grav$. To define the
product
$$\begin{CD}
\circ_i : \Grav\(m+1\)\o\Grav\(n+1\) @>>> \Grav\(m+n\) ,
\end{CD}$$
consider the embedding $j$ of
$\CM_{0,\{0,\dots,m\}}\times\CM_{0,\{0',\dots,n'\}}$ as a stratum of
$\Mbar_{0,\{0,\dots,\hat{\imath},\dots,m,1',\dots,n'\}}$, corresponding to
the joining of the point labelled $i$ in the curve
$\Sigma_1\in\CM_{0,\{0,\dots,m\}}$ to the point labelled $0'$ in the curve
$\Sigma_2\in\CM_{0,\{0',\dots,n'\}}$. Consider the Poincar\'e residue map
associated to this embedding \cite{Deligne}:
$$\begin{CD}
\Res : H^\bull(\CM_{0,\{0,\dots,\hat{\imath},\dots,m,1',\dots,n'\}}) @>>>
H^\bull(\CM_{0,\{0,\dots,m\}}\times\CM_{0,\{0',\dots,n'\}})
\end{CD}$$
Suitably suspending the adjoint of this map, we obtain the product
$\circ_i$ of $\Grav$. This construction makes it quite obvious that $\Grav$
satisfies the equivariance and associativity axioms of a cyclic operad.

Denote by $[x_1,\dots,x_n]$ the element of $\Grav(n)$ of degree $2-n$
corresponding to the standard basis vector of $H_0(\CM_{0,n+1})$. This
operation is graded antisymmetric, since $\SS_{n+1}$ acts on it by the sign
representation, and it is proved in \cite{weil} that this sequence of
operations generates $\Grav$, and that all relations are generated by the
quadratic relations of \eqref{gravity}. Note that in that paper, we work with
the operad $\Lambda^{-1}\Grav$, and the generators are all in degree
$1$: the relationship between these two sets of generators at the level of
algebras is
$$
\{x_1,\dots,x_n\} = (-1)^{(n-1)|x_1|+(n-2)|x_2|+\dots+|x_{n-1}|}
\Sigma^{-1} [\Sigma x_1,\dots,\Sigma x_n] .
$$
In that paper, the operad structure is constructed in a different, though
equivalent, way, using $\C^\times$-equivariant homology.

\subsection{$\Grav$ as a mixed Hodge operad}
Operads (and, more specifically, stable cyclic operads) may be defined in
any symmetric monoidal category with colimits. Up to this point, we have
concentrated on the examples of operads in the categories of chain
complexes (differential graded operads) and topological spaces (topological
operads). However, the category of mixed Hodge complexes \cite{BZ} is a
symmetric monoidal category with colimits, with graded tensor product as
the monoidal structure, and operads in this category are called mixed Hodge
operads. In fact, the mixed Hodge operads which most concern us are pure
and have vanishing differential.

The gravity operad is an example of a mixed Hodge operad. This
carries a unique mixed Hodge structure compatible with the Poincar\'e
residue maps which define the products in $\Grav$:
$$
\Grav\(n\) = \Sigma^{3-n} \sgn_n \o H_\bull(\CM_{0,n},\C(n-3)) .
$$
Here, $\C(n-3)$ is the Tate Hodge structure, which is a line with Hodge
numbers $(n-3,n-3)$.

\subsection{The hypercommutative operad}
Let $\mbar$ be the stable cyclic $\SS$-module
$$
\mbar\(n\) = \begin{cases} H_\bull(\Mbar_{0,n}) , & n\ge3 , \\ 0 , & n<3 .
\end{cases}$$
The cyclic $\SS$-space $\Mbar\(n\)=\Mbar_{0,n}$ is a topological cyclic
operad: the product is given by gluing stable curves together at marked
points. It follows that $\mbar$ is a cyclic operad. Kontsevich and Manin
found \cite{KM} that algebras over $\mbar$ are just hypercommutative
algebras in the sense of \eqref{hypercommutative}, where the operation
$(x_1,\dots,x_n)\in\mbar\(n+1\)$ corresponds to the fundamental class
$[\Mbar_{0,n+1}]\in H_{2(n-2)}(\Mbar_{0,n+1})$. The fact that the
operations $(x_1,\dots,x_n)$ generate $\mbar$ is an elegant restatement of
Theorem \ref{Keel:1}. In Proposition \ref{orthogonal}, we will give a new
proof (obtained jointly with Kontsevich) of the relations between the
generators of $\mbar$, which relies on the duality between hypercommutative
algebras and gravity algebras, together with our explicit presentation
\eqref{gravity} of the gravity operad $\Grav$.

\subsection{A spectral sequence of Deligne}
Let $M$ be a smooth projective variety of complex dimension $n$, and let
$\{D_1,\dots,D_N\}$ be a sequence of smooth divisors with normal crossings;
we denote their union by $D$. The sheaf of logarithmic differential forms
$\E^\bull_M(\log D)$ on $M$ is generated over the sheaf of differential
forms $\E^\bull_M$ by symbols $d(\log f)$, where $f$ is a section of
$\CO(D)$, subject to the relations
$$
d(\log fg) = d(\log f) + d(\log g) \quad\text{and}\quad
f\*d(\log f) = df.
$$
The sheaf $\E^\bull_M$ has a differential $d$, characterized by
$d(d(\log f))=0$, and we have the fundamental isomorphism
$$
H^\bull(U) \cong H^\bull(M,\E^\bull_M(\log D)) ,
$$
where $U=M\setminus D$.

The sheaf $\E^\bull_M(\log D)$ is filtered by subsheaves
$$
W_k\E^\bull_M(\log D) = \operatorname{span}_\CO\{
d(\log f_1)\.\dots\.d(\log f_i) \mid i\le k \} .
$$
Let $j^k:D^k\hookrightarrow M$ be the embedding of the closed subvariety
$$
D^k = \coprod_{i_1<\dots<i_k} D_{i_1} \cap \dots \cap D_{i_k} ,
$$
and let $\epsilon_k$ be the locally constant line bundle over $D^k$, which
over the component $D_{i_1}\cap\ldots\cap D_{i_k}$ equals the determinant
line $\det(\{i_1,\dots,i_k\})$. There is a canonical quasi-isomorphism
$$
\gr^W_k \E^\bull_M(\log D) \simeq
\Sigma^{-k}j^k_*\E^\bull_{D^k}\o\epsilon_k .
$$
The associated spectral sequence
$$
E_1^{-p,q} = H^{2p+q}(D^p,\epsilon_p) \abuts
E_\infty^{-p,q} = \gr^W_pH^{-p+q}(U)
$$
carries a Hodge filtration $F$, induced by the Hodge filtration of
$\E^\bull_M(\log D)$, and by the principal of two types, $E_2=E_\infty$
(\cite{Deligne}, Section 3.2). The weight filtration induced on $H^n(U)$ by
$W$ defines, up to translation, its mixed Hodge structure:
$\gr^W_pH^{-p+q}(M)$ carries a pure Hodge structure of weight $q$.

If a finite group $\Gamma$ acts on $M$, preserving $U\subset M$, this
spectral sequence carries an action of $\Gamma$ compatible with its action
on $H^\bull(U)$.

\subsection{The cohomology ring of $\CM_{0,n}$}
We now describe the cohomology ring of the moduli space $\CM_{0,n}$. Our
main tool is Arnold's description of the cohomology ring of the
configuration space $\C^n_0$ \cite{Arnold}.

If $1\le j\ne k\le n$, let $\om_{jk}$ be the logarithmic differential form
on the configuration space $\C^n_0$ given by the formula
$$
\om_{jk} = \frac{d\log(z_j-z_k)}{2\pi i} .
$$
Note that the cohomology class of $\om_{jk}$ is integral.
\begin{proposition} \label{Arnold}
The cohomology ring $H^\bull(\C^n_0,\Z)$ is the graded commutative ring
with generators $[\om_{jk}]$, and relations $\om_{jk}=\om_{kj}$ and
$\om_{ij}\om_{jk}+\om_{jk}\om_{ki}+\om_{ki}\om_{ij}=0$. The symmetric group
$\SS_n$ acts on $H^\bull(\C^n_0,\Z)$ through its action on the generators
$\sigma\*\om_{ij}=\om_{\sigma(i)\sigma(j)}$.
\end{proposition}
\begin{pf}
The Serre spectral sequence for the fibration
$$\begin{CD}
\C\setminus\{z_1,\dots,z_n\} @>>> \C^{n+1}_0 @>>> \C^n_0
\end{CD}$$
defined by projecting $(z_1,\dots,z_{n+1})$ to $(z_1,\dots,z_n)$ collapses
at $E_2$, and the monodromy of $\pi_1(\C^n_0)$ on
$H^\bull(\C\setminus\{z_1,\dots,z_n\})$ is trivial. The proof now
proceeds by induction on $n$.
\end{pf}

\begin{corollary} \label{equivariant}
The cohomology ring $H^\bull(\CM_{0,n+1},\C)$ may be identified with the
kernel of the differential $\iota$ on $H^\bull(\C^n_0,\C)$ whose action
on the generators is $\iota\om_{jk}=1$.
\end{corollary}
\begin{pf}
The isotropy group of the point $\infty\in\CP^1$ under the action of
$\PSL(2,\C)$ is
$$\textstyle
\Aff(\C) = \Bigl\{ \left( \begin{smallmatrix} a & b \\ 0 & a^{-1}
\end{smallmatrix} \right) \mathop{\Big|}
a\in\C^\times , b\in\C \Bigr\} \subset \PSL(2,\C) .
$$
Since $\PSL(2,\C)$ acts transitively on $\CP^1$, we see that
$\CM_{0,n+1}\cong\C^n_0/\Aff(\C)$. But the group $\Aff(\C)$ is homotopy
equivalent to the circle group
$$\textstyle
\Bigl\{ \left( \begin{smallmatrix} a & 0 \\ 0 & a^{-1} \end{smallmatrix}
\right) \mathop{\Big|} |a|=1 \Bigr\} ,
$$
giving a homotopy equivalence $\C^n_0\simeq\CM_{0,n+1}\times S^1$.
This allows us to identify the cohomology of $\CM_{0,n+1}$ with the
$S^1$-equivariant cohomology of $\C^n_0$.

The infinitesimal generator of the circle action on $\C^n_0$ is
the vector field
$$
T = 2\pi i \sum_{k=1}^n \bigl( z_k \p_k - \bar{z}_k \bar\p_k \bigr) ,
$$
whose contraction with a generator $\om_{jk}$ is $\om_{jk}(T)=1$.
\end{pf}

The above result leads to yet another realization of the $\SS_n$-module
structure on $\Lie\(n\)\cong\Lie(n-1)$:
$$
\Lie\(n\) \cong \sgn_n \o H_{n-3}(\CM_{0,n}) .
$$

\subsection{Application to $\CM_{0,n}\subset\Mbar_{0,n}$} \label{Deligne}
We now apply Deligne's spectral sequence with $U=\CM_{0,n}$ and
$M=\Mbar_{0,n}$. (This spectral sequence is also discussed in \cite{GK},
Section 3.4.5.) Denote the closure of the open stratum $\CM\(T\)$ by
$\Mbar\(T\)$. Then the divisors are the closed strata $\CM\(T\)$,
$T\in\CT_1\(n\)$, while $D^p$ is the union of the closed strata
$\Mbar\(T\)$, $T\in\CT_p\(n\)$. The restriction of $\epsilon_p$ to
$\Mbar\(T\)$ equals $\det(\Edge(T))$; as this is naturally isomorphic to
$\det(T)$, we see that
$$
E_1^{-p,q} \cong \bigoplus_{T\in\CT_p\(n\)}
H^{-2p+q}(\Mbar\(T\),\det(T)) .
$$
The differential $\begin{CD}d_1:E_1^{-p,q}@>>>E_1^{-p+1,q}
\end{CD}$ is easy to describe: it
is the composition
$$\begin{CD}
\displaystyle
\bigoplus_{T\in\CT_p\(n\)} H^{-2p+q}(\Mbar\(T\),\det(T)) @.
\displaystyle
\bigoplus_{T\in\CT_{p-1}\(n\)} H^{-2p+2+q}(\Mbar\(T\),\det(T)) \\
@| @| \\
\displaystyle
\bigoplus_{T\in\CT_p\(n\)} H_{2(n-3)-q}(\Mbar\(T\),\det(T)) @>>>
\displaystyle
\bigoplus_{T\in\CT_{p-1}\(n\)} H_{2(n-3)-q}(\Mbar\(T\),\det(T))
\end{CD}$$
where the vertical isomorphisms are induced by Poincar\'e duality, and the
bottom arrow is the map induced on the homology groups by the inclusion of
the $p$-codimensional closed strata of $\Mbar_{0,n}$ into the
$(p-1)$-codimensional closed strata.

The key to unlocking this spectral sequence is the following lemma, which
shows that
$$
E_2^{-k,2k} \cong H^k(\CM_{0,n}) ,
$$
while $E_2^{pq}=0$ if $2p+q\ne0$.

\begin{lemma} \label{mixed}
The mixed Hodge structure of $H^k(\CM_{0,n})$, $n\ge3$, is pure of weight
$2k$.

In degree $i$, $\Grav\(n\)$ equals the $\SS_n$-module $\sgn_n\o
H_{i+n-3}(\CM_{0,n})$, with pure Hodge structure of weight $-2i$.
\end{lemma}
\begin{pf}
By Proposition \ref{Arnold}, the cohomology ring of $\C^n_0$ is generated
by the logarithmic differential forms $\om_{ij}$; it follows that the mixed
Hodge structure of $H^k(\C^n_0)$ is pure of weight $2k$. By Corollary
\ref{equivariant}, there is an injection of the cohomology of $\CM_{0,n+1}$
into the cohomology of $\C^n_0$, induced by the quotient map
$\begin{CD}\C^n_0@>>>\CM_{0,n+1}
\end{CD}$, and the result follows.
\end{pf}

Thus, for $q$ even, the $E_1$-term of the spectral sequence gives rise to a
resolution of the graded vector space $H^\bull(\CM_{0,n})$
\begin{equation} \label{bar:even} \begin{CD}
0 @>>> H^p(\CM_{0,n}) @>>>
\bigoplus_{T\in\CT_p\(n\)} H^0(\Mbar\(T\),\det(T)) @>>>
\bigoplus_{T\in\CT_{p-1}\(n\)} H^2(\Mbar\(T\),\det(T)) @>>> \dots
\end{CD} \end{equation}
For $q$ odd, we obtain the exact sequence
$$\begin{CD}
\dots @>>> \bigoplus_{T\in\CT_1\(n\)} H^{q-2}(\Mbar\(T\),\det(T)) @>>>
H^q(\Mbar_{0,n}) = \bigoplus_{T\in\CT_0\(n\)}
H^q(\Mbar\(T\),\det(T)) @>>> 0 .
\end{CD}$$
which shows, by induction on odd $q$, that the cohomology groups
$H^q(\Mbar_{0,n})$ vanish if $q$ is odd.

\subsection{Proof of Theorem \ref{Keel:1}} \label{Keel:proof}
It is clear that the result holds for $n=3$, since $\Mbar_{0,3}$ is a
point. We now argue by induction on $n$. As a consequence of the
surjectivity of the differential in the exact sequence \eqref{bar:even}
$$\begin{CD}
\bigoplus_{T\in\CT_1\(n\)} H^{2p-2}(\Mbar\(T\),\det(T)) @>>>
H^{2p}(\Mbar_{0,n}) = \bigoplus_{T\in\CT_0\(n\)}
H^{2p}(\Mbar\(T\),\det(T))
\end{CD}$$
for $p\le n-2$, we see that all homology classes of $\Mbar_{0,n}$ except
the fundamental class are supported on the closures of strata of
codimension $1$. Such a closed stratum is isomorphic to
$\Mbar_{0,i}\times\Mbar_{0,j}$ where $i+j=n+2$ and $i,j\ge3$, allowing us
to apply the induction. \qed

\subsection{The dimension of $H^2(\Mbar_{0,n})$}
Another consequence of \eqref{bar:even} is a simple formula for the
dimension $\dim H^2(\Mbar_{0,n})$ of the Picard variety of $\Mbar_{0,n}$:
\begin{eqnarray*}
\dim H^2(\Mbar_{0,n}) &=& \sum_{k=3}^{n-1} \binom{n-1}{k}
= 2^{n-1} - \frac{n^2-n+2}{2} \\
&=& \binom{n}{n-4} + \binom{n}{n-6} + \dots
\end{eqnarray*}
\begin{pf}
When $p=1$, the short exact sequence \eqref{bar:even} becomes
$$\begin{CD}
0 @>>> H^1(\CM_{0,n}) @>>>
\bigoplus_{T\in\CT_1\(n\)} H^0(\Mbar\(T\),\det(T)) @>>>
H^2(\Mbar_{0,n}) @>>> 0 .
\end{CD}$$
By Corollary \ref{equivariant}, the dimension of $H^1(\CM_{0,n})$ equals
the coefficient of $-t$ in $(1-2t)\dots(1-(n-2)t)$, or
$$
2 + \dots + (n-1) = \binom{n-1}{2} - 1 .
$$
(See \eqref{poincare-moduli} for more details of this calculation: we will
actually show that $H^1(\CM_{0,n})$ is isomorphic to the irreducible
$\SS_n$-module $V_{n-2,2}$.) Each tree in $\CT_1\(n\)$ contributes a copy
of $\C$ to $\bigoplus_{T\in\CT_1\(n\)} H_{2(n-4)}(\Mbar\(T\))$. Let $T$ be
such a tree and consider the set $S$ of the external edges attached to one
vertex of $T$. We see that trees with two vertices correspond to subsets
$S\subset\{1,\dots,n\}$ where $2\le|S|\le n-2$, where we identify the trees
corresponding to the subsets $S$ and $S^c$, the complement of $S$. Thus,
$$
|\CT_1\(n\)| = \frac12 \sum_{k=2}^{n-2} \binom{n}{k}
= \frac12 (2^n-2n-2) = 2^{n-1} - n - 1 .
$$
The result follows easily.
\end{pf}

This formula may be compared to the dimension of $H^2(\Mbar_{g,n})$, $g>2$,
which follows from the work of Arbarello and Cornalba \cite{AC}:
$$
\dim H^2(\Mbar_{g,n}) = 2^{n-1} (g+1) + n + 1 .
$$
We see that this formula is correct for $g=0$ up to a polynomial error.

The above dimension formula may be refined, using the realization of
$\Mbar_{0,n}$ as an iterated blowup, to show that as an $\SS_n$-module,
$H^2(\Mbar_{0,n})$ is the direct sum of the suitable exterior powers of the
permutation representation $\C^n$ of $\SS_n$.

Note that there are the same number of $n$-linear relations \eqref{gravity}
among the brackets $[x_1,\dots,x_k]$ generating the gravity operad as $\dim
H^2(\Mbar_{0,n+1})$; as we will see, this is no coincidence.

\section{Koszul operads}

In this section, we prove our main theorem, the duality of the
hypercommutative and gravity operads. To do this, we must generalize
Ginzburg and Kapranov's theory of Koszul operads \cite{GK} so that it
applies to operads which are not necessarily generated by bilinear
operations. First, we recall their cobar construction for operads, an
analogue of Hochschild's bar construction for associative algebras.

The dual of an operad is only defined up to homotopy, and is represented by
the cobar operad. However, there is a class of operads, the Koszul operads,
for which there is a particularly nice dual, whose generators are in
one-to-one correspondence with those of the original operad. A Koszul
operad is quadratic (the relations among its generators are bilinear), as
is its dual, and the relations in the dual operad may be characterized as
the orthogonal complement of those of the original operad.

\subsection{Free operads and trees}
We now recall from \cite{cyclic} the structure of the free cyclic operad
$\TT_+\v$ generated by a cyclic $\SS$-module $\v$. There is an analogous
construction for operads, for which we refer to \cite{n-algebras}. From now
on, all cyclic $\SS$-modules which we discuss will be stable.

If $\v$ is a (stable) cyclic $\SS$-module, let $\TT_+\v$ be the
(stable) cyclic $\SS$-module defined by
$$
\TT_+\v\(n\) = \bigoplus_{T\in\CT\(n\)} \v\(T\) ,
$$
where $\v\(T\)=\bigotimes_{v\in\VV(T)} \v\(\Leg(v)\)$. Note that
$\TT_+\v$ is graded by subspaces
$$
\TT_i\v\(n\) = \bigoplus_{T\in\CT_i\(n\)} \v\(T\) .
$$
Then $\TT_+$ is an endofunctor in the category of (stable) cyclic
$\SS$-modules.

There is a natural structure of a triple on the functor $\TT_+$:
\begin{enumerate}
\item since $\TT_+\TT_+$ is a sum over trees, each vertex of which is
itself a tree, the product of the triple is a natural transformation
from $\TT_+\TT_+$ to $\TT_+$ obtained by gluing the trees at the vertices
into the larger tree;
\item the unit of the triple is the natural transformation from the
identity functor to $\TT_+$ induced by the inclusion
$\CT_0\(n\)\subset\CT\(n\)$.
\end{enumerate}
The following theorem is a melding of results from \cite{cyclic} and
\cite{modular}.

\begin{theorem}
A (non-unital, stable) cyclic operad is the same thing as a $\TT_+$-algebra
in the category of (stable) cyclic $\SS$-modules.
\end{theorem}

\subsection{The cobar construction for operads}
The cobar construction $\BB$, introduced by Ginzburg and Kapranov
\cite{GK}, is a contravariant functor on the category of operads. We study
here the slight variant of this functor which acts on the category of
(non-unital, stable) cyclic operads.

The dual $V^*$ of a chain complex $V$ is defined as follows:
$\begin{CD}V^*_i=(V_{-i})^*$, and $\delta^*:V^*_i@>>>V^*_{i-1}
\end{CD}$ is the adjoint of
$\begin{CD}\delta:V_{-i+1}@>>>V_{-i}\end{CD}$.

If $\v$ is a stable cyclic $\SS$-module, denote by $\v^\Dual$ the
stable cyclic $\SS$-module
$$
\v^\Dual\(n\) = \Sigma^{n-3}\sgn_n\o\v\(n\)^* .
$$
This functor is an involution on the category of stable cyclic
$\SS$-modules, that is, $(\v^\Dual)^\Dual$ is naturally isomorphic to
$\v$.

The cobar operad $\BB\a$ of a (non-unital, stable) cyclic operad is
obtained by perturbing the differential of the free cyclic operad
$\TT_+\a^\Dual$ by a differential $\p$ which reflects the operad structure
of $\a$, and is defined as follows.

If $T\in\CT\(n\)$, and $e$ is an edge of $T$, denote by $T/e$ the tree in
which $e$ is contracted to a point: thus, $T/e$ has one fewer vertices, and
one fewer edges, than $T$. There is a natural map of degree $0$
$$\begin{CD}
\p_{T/e} : \a\(T\) @>>> \a\(T/e\)
\end{CD}$$
induced by composition in the operad $\a$ along the edge $e$. This induces
a map
$$\begin{CD}
\p^\Dual_{T/e} : \a^\Dual\(T/e\) @>>> \a^\Dual\(T\)
\end{CD}$$
of degree $-1$. We now define the differential $\p$ to be the operator
whose matrix element from $\a^\Dual\(\tilde{T}\)\subset\BB\a$ to
$\a^\Dual\(T\)\subset\BB\a$ is the sum of the operators $\p^\Dual_{T/e}$
over internal edges $e$ such that $T/e$ is isomorphic to $\tilde{T}$.

Paying careful attention to the signs coming from the suspensions, one
shows that the differential $\p$ satisfies the formulas
$\p^2=\delta\p+\p\delta=0$, and hence that $\delta+\p$ is a differential on
$\TT_+\a^\Dual$. It is also not hard to show that $\p$ is compatible with
the cyclic operad structure of $\TT_+\a^\Dual$, so that
$\BB\a=(\TT_+\a^\Dual,\delta+\p)$ is an operad.

The properties of the resulting functor are summarized by the following
theorem.
\begin{theorem} \label{cobar}
\textup{(1)} The cobar construction is a homotopy functor, that is, if
$\begin{CD}f:\a@>>>\b
\end{CD}$ is a homotopy equivalence, then so is $\begin{CD}
\BB f:\BB\a@>>>\BB\b\end{CD}$.

\noindent \textup{(2)} There is a natural transformation from $\BB\BB$ to the
identity functor, and the resulting map $\begin{CD}\BB\BB\a@>>>\a
\end{CD}$ is a
homotopy equivalence for all $\a$.
\end{theorem}
\begin{pf}
The homotopy invariance of $\BB$ is easy to see by a double complex
argument. The natural map from $\BB\BB\a$ to $\a$, which is projection onto
the summand $\a\cong(\a^\Dual)^\Dual\subset\TT_+(\TT_+\a^\Dual)^\Dual$, is
shown to be a homotopy equivalence of operads in Theorem 3.2.16 of
\cite{GK}.
\end{pf}

If $\a$ and $\b$ are operads and $\begin{CD}\Phi:\BB\a@>>>\b
\end{CD}$ is a morphism of
operads, there is a bar construction $\BC_\Phi$ on $\a$-algebras, defined
for an $\a$-algebra $A$ by twisting the differential on the chain complex
$$
A \oplus
\Sigma^{-1} \bigoplus_{n=3}^\infty \Hom_{\SS_n}(\b(n),(\Sigma A)^{\o n})
$$
in such a way as to reflect the $\a$-algebra structure of $A$. (See
\cite{n-algebras} for details.) When $\b=\BB\a$ and $\Phi$ is the identity
map, we denote the resulting functor $\BC$. Let $QA$ be the complex of
indecomposables, obtained by taking the cokernel of the map
$$\begin{CD}
\bigoplus_{n=3}^\infty \rho_n :
\bigoplus_{n=3}^\infty \a(n)\o_{\SS_n} A^{\o n} @>>> A .
\end{CD}$$
The following theorem is proved in Chapter 3 of \cite{n-algebras}.
\begin{theorem} \label{Bar}
\textup{(1)} There is a natural transformation of functors $\begin{CD}
\BC@>>>Q\end{CD}$, such
that if $A$ is a free algebra, the morphism $\begin{CD}\BC A@>>>QA
\end{CD}$ is a homotopy
equivalence.

\noindent \textup{(2)} The functors $\BC_\Phi$ are homotopy functors: if
$\begin{CD}f:A@>>>B
\end{CD}$ is a homotopy equivalence, then so is $\begin{CD}\BC_\Phi f:\BC_\Phi
A@>>>\BC_\Phi B\end{CD}$.
\end{theorem}

If $\begin{CD}\Phi:\BB\a@>>>\b
\end{CD}$ is a homotopy equivalence of operads, the natural
morphism $\begin{CD}\BC A@>>>\BC_\Phi A
\end{CD}$ is a homotopy equivalence for all
$\a$-algebras $A$. Thus, the functor $\BC_\Phi$ is a left derived functor
$\L Q$ of the indecomposable functor $Q$, that is, a homotopy functor
homotopy equivalent to $Q$ on free $\a$-algebras. It is proved in
\cite{GK} that there are natural homotopy equivalences of operads
$$\begin{CD}
\BB\Ass @>>> \Ass \quad \BB\Com @>>> \Lie \quad \BB\Lie @>>> \Com .
\end{CD}$$
The bar construction associated to the first of these homotopy equivalences
is, up to a shift in degree, the Hochschild bar construction on associative
algebras, while the bar constructions associated to the other two homotopy
equivalences are the functors $\L Q$ on commutative and Lie algebras
discussed in the introduction, equal, up to a shift in degree, to the
Harrison and Chevalley-Eilenberg complexes respectively. Thus, the duality
result Theorem \ref{Harrison} is seen to be a special case of Theorem
\ref{Bar}.

In Theorem \ref{Kontsevich}, we will prove that $\mbar$-algebras are the
same thing as hypercommutative algebras in the sense of
\eqref{hypercommutative}. Thus the duality between hypercommutative and
gravity algebras announced in the introduction follows from Theorem
\ref{Bar} combined with the following generalization of Theorem 4.25 of
\cite{GK}.
\begin{theorem} \label{duality}
There is a natural homotopy equivalence of operads $\begin{CD}
\BB\mbar@>>>\Grav\end{CD}$.
\end{theorem}
\begin{pf}
Let $\v$ be the $\SS$-module obtained by summing the short exact sequences
\eqref{bar:even} (minus the terms $H^p(\CM_{0,n})$), placing the summand
$$
\bigoplus_{T\in\CT_p\(n\)} H^q(\Mbar\(T\),\det(T))
$$
of $\v\(n\)$ in degree $2(n-3)-p-q$. Using Poincar\'e duality, we see that
$$
\v\(n\) \cong \bigoplus_{p=0}^{n-3}
\bigoplus_{T\in\CT_p\(n\)} \Sigma^p \det(T) \o \mbar\(T\) .
$$
Furthermore, there is a natural homotopy equivalence
$\begin{CD}\Sigma^{2(n-3)}H^\bull(\Mbar_{0,n})@>>>\v\(n\)
\end{CD}$, which induces a homotopy
equivalence $\begin{CD}\v^\Dual@>>>\Grav\end{CD}$.

By the isomorphism \eqref{det}, the $\SS$-module $\v^\Dual$ may be rewritten
as
\begin{eqnarray*}
\v^\Dual\(n\) &\cong& \bigoplus_{p=0}^{n-3} \Sigma^{n-3-p}
\bigoplus_{T\in\CT_p\(n\)} \bigotimes_{v\in\VV(T)} \det(\Leg(v)) \o
\mbar^*\(\Leg(v)\) \\
&\cong& \bigoplus_{p=0}^{n-3} \Sigma^{n-3-p}
\bigoplus_{T\in\CT_p\(n\)} \bigotimes_{v\in\VV(T)} \Sigma^{3-|v|}
\mbar^\Dual\(\Leg(v)\) \\
&\cong& \bigoplus_{p=0}^{n-3} \Sigma^{n-3-p+\sum_{v\in\VV(T)}(3-|v|)}
\bigoplus_{T\in\CT_p\(n\)} \mbar^\Dual\(T\) .
\end{eqnarray*}
If $T\in\CT_p\(n\)$, we have
$$
\sum_{v\in\VV(T)} (3-|v|) = 3 - n + p ,
$$
showing that
$$
\v^\Dual\(n\) \cong \bigoplus_{p=0}^{n-3} \bigoplus_{T\in\CT_p\(n\)}
\mbar^\Dual\(T\) \cong \BB\mbar\(n\) .
$$
A little diagram chasing shows that the differentials of the $\SS$-modules
$\v^\Dual$ and $\BB\mbar$ are the same, and that the resulting homotopy
equivalence between $\Grav$ and $\BB\mbar$ is compatible with the operad
structures.
\end{pf}

\subsection{The cobar construction for mixed Hodge operads}
The free operad functor $\v\mapsto\TT\v$ and the functor $\begin{CD}
\v@>>>\v^\Dual\end{CD}$
have analogues in the category of mixed Hodge $\SS$-modules, defined in
precisely the same way as in the category of $\SS$-modules. (We recall that
the dual $\v^*$ in the category of mixed Hodge complexes reverses the
weight filtration, sending complexes of weight $k$ to complexes of weight
$-k$.) This allows us to extend the cobar construction to the category of
mixed Hodge operads, by the same definition as in the category of dg
operads.

If we follow through the proof of Theorem \ref{duality} paying attention to
the mixed Hodge structures, we see that the homotopy equivalence
$\begin{CD}\BB\mbar@>>>\Grav
\end{CD}$ is compatible with the Hodge structures of $\Grav$ and
$\mbar$, where $\mbar$ carries the natural (pure) Hodge structure coming
from its realization as the cohomology of the smooth K\"ahler manifold
$\Mbar_{0,n}$. This observation will be essential in our calculation of the
defining relations of the operad $\mbar$.

\subsection{Quadratic operads}
We now generalize Ginzburg and Kapranov's notion of a Koszul operad to
operads whose generators are not necessarily bilinear operations. Once
more, we restrict attention to stable cyclic operads.

An ideal $\b\subset\a$ of a cyclic operad is a cyclic $\SS$-submodule such
that for all operations $\circ_i$, $a\circ_ib$ is in $\b$ if either $a$ or
$b$ is. The intersection of two ideals is obviously an ideal. An ideal is
generated by a cyclic $\SS$-submodule $\r\subset\b$ if $\b$ is the
intersection of all ideals of $\a$ containing $\r$.

Let $\a$ be an operad, generated by a cyclic $\SS$-submodule $\v$. The pair
$(\a,\v)$ is a cyclic quadratic operad if the ideal $\begin{CD}
\ker(\TT\v@>>>\a)\end{CD}$ in
the free cyclic operad $\TT_+\v$ is generated by
$$\begin{CD}
\ker \Bigl( \TT_1\v=\bigoplus_{T\in\CT_1} \v\(T\) @>>> \a \Bigr) .
\end{CD}$$
The word quadratic is used here because $\v\(T\)$ is quadratic in $\v$ if
$T$ has one internal edge and hence two vertices. Thus, $\r$ is itself
quadratic in $\v$.

The cyclic operads $\a=\Ass$, $\Com$ and $\Lie$ are all quadratic, with
generating cyclic submodule $\v$, where
$$
\v\(n\) = \begin{cases} \a\(3\) , & n=3, \\ 0 , & n\ne3 . \end{cases}
$$
For example, $\Lie\(3\)$ is one-dimensional, spanned by $[a_1,a_2]$, and
the cyclic $\SS$-module of relations $\r$ is given by the formula
$$
\r\(n\) =
\begin{cases} \span\{ [a_1,[a_2,a_3]] , [a_2,[a_3,a_1]] \} , & n=4 , \\
0 , & n\ne4 . \end{cases}
$$

\subsection{The naive dual of a quadratic operad}
If $\a$ is a cyclic quadratic operad, the naive dual $\a^!$ of $\a$ is the
cokernel of the composition
$$\begin{CD}
\psi : \BB\a @>\p>> \BB\a @>>> \TT_+\v^\Dual ,
\end{CD}$$
where the second arrow is the surjection of cyclic $\SS$-modules
$\begin{CD}\BB\a@>>>\TT_+\v^\Dual
\end{CD}$ induced by the inclusion of cyclic $\SS$-modules
$\v\subset\a$.

\begin{definition}
If $\a$ is a cyclic quadratic operad, there is a natural morphism of
operads $\begin{CD}\Phi:\BB\a@>>>\a^!\end{CD}$, induced by the surjection
$\begin{CD}\BB\a@>>>\TT_+\v^\Dual
\end{CD}$. The operad $\a$ is Koszul if the surjection of
operads $\begin{CD}\BB\a@>>>\a^!
\end{CD}$ is a homotopy equivalence, or equivalently, if
$\begin{CD}\BB\a^!@>>>\a\end{CD}$ is.
\end{definition}

\begin{proposition} \label{orthogonal}
If $\a$ is a cyclic quadratic operad with generators $\v$ and relations
$\r$, let $\r^\perp$ be the kernel of the natural map from
$\TT_1\v^\Dual\cong(\TT_1\v^\Dual)^\Dual$ to $\r^\Dual$. Then $\a^!$ is a
cyclic quadratic operad with generators $\v^\Dual$ and relations $\r^\perp$.
\end{proposition}
\begin{pf}
It suffices to show that the image of $\BB\a$ in $\TT_+\v^\Dual$ under the
above composition is the ideal generated by $\r^\perp$. Denote by
$\a_k\subset\a$ the image of $\TT_k\v$ in $\a$ under the quotient map
$\begin{CD}\TT_+\v@>>>\a
\end{CD}$. Thus, $\a_0=\v$ and $\a_1=\TT_1\v/\r$. Observe that
$\a_1^\Dual\cong\r^\perp$.

If $T\in\CT\(n\)$, the summand $\a^\Dual\(T\)$ of $\BB\a\(n\)$ may be
thought of as the vector space spanned by decorations of the tree $T$, in
which each vertex of $T$ is assigned an element of $\a^\Dual$ of
appropriate valence. The map $\begin{CD}\psi:\BB\a@>>>\TT_+\v^\Dual
\end{CD}$ vanishes on such
a decorated tree unless the vertex decorations lie in $\v$ at all but one
vertex $v$, which is decorated by $a\in\a_1^\Dual$. The map $\psi$ applied
to this decorated tree produces a new decorated tree in which the vertex
$v$ is replaced by the tree underlying $a$ (and thus having one additional
edge). Thus, the image of $\psi$ is the ideal generated by $\r^\perp$.
\end{pf}

\begin{corollary}
A cyclic quadratic operad $\a$ is Koszul if and only if $\a^!$ is, and
$(\a^!)^!\cong\a$.
\end{corollary}

As examples of naive duals, we have $\Ass^!\cong\Ass$, $\Com^!\cong\Lie$
and $\Lie^!\cong\Com$. It is proved in \cite{GK} that the operads $\Ass$,
$\Com$ and $\Lie$ are Koszul. A non-cyclic example of a Koszul operad is
$\Braid$ \cite{n-algebras}, which satisfies
$\Braid^!\cong\Lambda^{-1}\Braid$.

The proof of the following theorem occupies the remainder of this
section. This theorem is joint work of the author and M. Kontsevich.
\begin{theorem} \label{Kontsevich}
Let $\v\subset\mbar$ be the cyclic $\SS$-submodule spanned by the
fundamental classes
$$
[\Mbar_{0,n}] \in H_{2(n-3)}(\Mbar_{0,n}) \subset \mbar\(n\) .
$$
The operad $\mbar$ is Koszul, with generators $\v$, and
$\mbar^!\cong\Grav$.
\end{theorem}
\begin{pf}
This theorem is proved using the duality between the mixed Hodge operads
$\Grav$ and $\mbar$, and the fact that the operad $\Grav$ is quadratic.

If $\a$ is a mixed Hodge operad, the natural homotopy equivalence
$\begin{CD}\BB\BB\a@>>>\a
\end{CD}$ of Theorem \ref{cobar} is a morphism of mixed Hodge
operads. It follows that we have a diagram in the category of mixed Hodge
operads
$$\begin{CD}
\BB\BB\mbar @>>> \BB\Grav \\
@VVV @. \\
\mbar @.
\end{CD}$$
in which both arrows are homotopy equivalences. The homology of the weight
$-2p$ summand of the complex $\BB\Grav\(n\)$, which by this argument is
isomorphic to $H_{2p}(\Mbar_{0,n})$, must be concentrated in degree $2p$;
from this, we see that this subcomplex is exact except at the last term,
giving a long exact sequence
\begin{equation} \label{bar:Hodge}\begin{CD}
\dots @>>> \bigoplus_{T\in\CT_{n-2-p}\(n\)} \Grav^\Dual\(T\)_{2p+1}
@>>> \bigoplus_{T\in\CT_{n-3-p}\(n\)} \Grav^\Dual\(T\)_{2p}
@>>> \mbar\(n\)_{2p} @>>> 0 .
\end{CD}\end{equation}

Let $\v^\Dual\subset\Grav$ be the cyclic $\SS$-module spanned by the generators
of $\Grav$; for each $n\ge3$, there is one generator, of degree $3-n$ and
weight $2(n-3)$, in $\Grav\(n\)$. Thus $\v\(n\)$ is spanned by an
element of degree $2(n-3)$ and weight $2(3-n)$. Taking $p=n-3$ in the long
exact sequence \eqref{bar:Hodge}, we see that $\v\(n\)$ may be identified
with $H_{2(n-3)}(\Mbar_{0,n})$.

Let $\r^\Dual\subset\Grav$ be the cyclic $\SS$-module spanned by the
elements of degree $4-n$ in $\Grav$. An $\SS_n$-module, $\r^\Dual\(n\)\cong
H_1\(\CM_{0,n}\)\o\sgn_n$; it is concentrated in degree $4-n$, has weight
$2(n-4)$, and dimension $\binom{n-1}{2}-1$. The case $p=n-4$ of the long
exact sequence \eqref{bar:Hodge} is the short exact sequence
$$\begin{CD}
0 @>>> \r\(n\) @>>> \bigoplus_{T\in\CT_1\(n\)} \v\(T\)
@>>> H_{2(n-4)}(\Mbar_{0,n}) @>>> 0 ,
\end{CD}$$
showing that the $\SS$-module $\r$ is a subset of the set of
relations for $\mbar$. Furthermore, there are no further relations, as may
be seen from \eqref{bar:Hodge}: for general $p\le n-3$, there is an exact
sequence
\begin{multline*}\begin{CD}
\dots @>>> \bigoplus_{T\in\CT_{n-2-p}\(n\)} \bigoplus_{v\in\VV(T)}
\r\(\Leg(v)\) \o \bigotimes_{w\in\VV(T)\setminus\{v\}} \v\(\Leg(w)\) @>>> \\
\bigoplus_{T\in\CT_{n-3-p}\(n\)} \v\(T\) @>>> H_{2p}\(n\)_{2p} @>>> 0 .
\end{CD}\end{multline*}
This shows that the operad $\mbar$ is quadratic, with generators $\v$ and
relations $\r$.

It remains to identify the $\SS$-module $\r$ with the set of relations
\eqref{hypercommutative} which hold in a hypercommutative algebra. This is
done in two parts: we first show that these relations are are in the
orthogonal complement of the relations which define a gravity algebra, and
thus form a subset of $\r$, and then show that they form a subspace of
$\r\(n\)$ of dimension at least $\binom{n-1}{2}-1$. Since $\r\(n\)$ itself
has dimension $\binom{n-1}{2}-1$, this completes the proof.

It is simple to check that the relations \eqref{hypercommutative} are
orthogonal to those which hold in a gravity algebra. Consider the relation
$$
G_0 = \pm [[a_i,a_j],a_1,\dots,\widehat{a_i},\dots,\widehat{a_j},\dots,a_k] .
$$
The inner product of this relation with the relation
$$
H = \sum_{S_1\coprod S_2=\{1,\dots,n\}} \pm ((a,b,x_{S_1}),c,x_{S_2})
- \sum_{S_1\coprod S_2=\{1,\dots,n\}} \pm (a,(b,c,x_{S_1}),x_{S_2})
$$
vanishes, since only the terms with $S_1=\emptyset$ can contribute: they
each contribute a term $1$, but with opposite sign.

Turning now to the relation
$$
G_\ell = \pm
[[a_i,a_j],a_1,\dots,\widehat{a_i},\dots,\widehat{a_j},\dots,a_k]
- [[a_1,\dots,a_k],b_1,\dots,b_\ell] , \quad \ell>0 ,
$$
we see that there are three cases to consider:
\begin{enumerate}
\item none of the letters $a,b,c$ lie in the set $\{b_1,\dots,b_\ell\}$, in
which case the inner product of relation $H$ with the above relation again
vanishes, for the same reason as when $\ell=0$;
\item one of the letters, say $c$, lies in the set $\{b_1,\dots,b_\ell\}$,
in which case the only terms having a non-zero inner product with $H$ are
$[[a,b],c,x_S]$ and $[[a_1,\dots,a_k],b_1,\dots,b_\ell]$, whose
contributions, each equal to $1$, cancel;
\item two or three of the letters $a,b,c$ lie in the set
$\{b_1,\dots,b_\ell\}$, in which case the inner product of the above
relation with each term of $H$ vanishes.
\end{enumerate}

Finally, we check that the space of relations among $n$ letters in a
hypercommutative algebra has dimension at least $\binom{n}{2}-1$. Consider
the projection of these relations into the space $\a(n)$ of all quadratic
words in the generators of the hypercommutative operad of the form
$$
((x_i,x_j),x_1,\dots,\widehat{x_i},\dots,\widehat{x_j},\dots,x_n) .
$$
(Note that this subspace of $\TT\v^\Dual\(n+1\)$ is not
$\SS_{n+1}$-invariant, but only $\SS_n$-invariant.) The dimension of
$\a(n)$ is $\binom{n}{2}$, and the relations \eqref{hypercommutative} project
in $\a(n)$ to relations
$$
((a,b),c,x_S) = ((a,c),b,x_S) .
$$
Clearly, the quotient of $\a(n)$ by these relations is
one-dimensional. This completes the proof that \eqref{hypercommutative} are
all the relations in the operad $\mbar$.
\end{pf}

\section{The equivariant Poincar\'e polynomials of $\CM_{0,n}$ and
$\Mbar_{0,n}$}

In this section, we use the results of \cite{modular} to calculate the
character of the $\SS_n$-modules $H_i(\CM_{0,n})$ and
$H_i(\Mbar_{0,n})$. By and large, the results of this section are
independent of the rest of this paper.

\subsection{Symmetric functions}
Let $\Lambda$ be the ring of symmetric functions: this is the limit
$$
\Lambda = \varprojlim \Z\[x_1,\dots,x_k\]^{\SS_k} .
$$
Then $\Lambda$ is the ring $\Z\[h_1,h_2,\dots\]$ of power series in the
complete symmetric functions
$$
h_n(x_i) = \sum_{i_1\le\dots\le i_n} x_{i_1}\dots x_{i_n} ,
$$
and $\Lambda_\Q=\Lambda\o\Q$ is a power series ring $\Q\[p_1,p_2,\dots\]$
in the power sums
$$
p_n(x_i) = \sum_i x_i^n .
$$

If $\sigma\in\SS_n$ has cycles of length
$\lambda_1\ge\dots\ge\lambda_\ell$, its cycle index $\psi(\sigma)$ is the
monomial $p_{\lambda_1}\dots p_{\lambda_\ell}$. If $V$ is
an $\SS_n$-module, its characteristic is the symmetric function
$$
\ch_n(V) = \frac{1}{n!} \sum_{\sigma\in\SS_n} \Tr_V(\sigma) \psi(\sigma) .
$$
It may be shown that $\ch_n(V)\in\Lambda\subset\Lambda_\Q$, and that the
characteristics of the irreducible representations of the symmetric groups
$\SS_n$, $n\ge0$, form a basis of $\Lambda$ over $\Z$, called the Schur
functions \cite{Macdonald}. For example, $h_n$ is the characteristic of the
trivial representation of $\SS_n$.

Define the Poincar\'e characteristic of an $\SS$-module to be
$$
\ch_t(\v) = \sum_{n=0}^\infty \sum_{i=0}^\infty (-t)^i \, \ch_n(\v_i(n))
\in \Lambda\(t\) .
$$
Setting $t=1$, we obtain the (Euler-Frobenius) characteristic $\ch(\v)$.
For example,
$$
\ch(\Com^+) = \sum_{n=1}^\infty h_n = \exp\Bigl( \sum_{n=1}^\infty
\frac{p_n}{n} \Bigr) - 1 , \quad\text{and}\quad
\ch(\Ass^+) = \sum_{n=1}^\infty p_1^n = \frac{p_1}{1-p_1} .
$$

\subsection{Plethysm}
Consider the ring $\Lambda\(t\)$ of power series in a variable $t$ with
coefficients in $\Lambda$. There is an associative product on
$\Lambda\(t\)$, called plethysm and denoted $f\circ g$, characterized by
the formulas
\begin{enumerate}
\item $(f_1+f_2)\circ g=f_1\circ g+f_2\circ g$;
\item $(f_1f_2)\circ g=(f_1\circ g)(f_2\circ g)$;
\item if $f=f(t,p_1,p_2,\dots)$, then $p_n\circ f=f(t^n,p_n,p_{2n},\dots)$,
and $t\circ f=t$.
\end{enumerate}
The following formula generalizes its analogue for ungraded $\SS$-modules,
proved in \cite{Macdonald}:
$$
\ch_t(\v\circ\w) = \ch_t(\v)\circ\ch_t(\w) .
$$

The operation
$$
\Exp(f) = \sum_{n=0}^\infty h_n\circ f
$$
plays the role for symmetric functions that exponentiation does for
power series. The inverse of $\Exp$ is the operation
$$
\Log(f) = \sum_{n=1}^\infty \frac{\mu(n)}{n} \log(p_n\circ f) ,
$$
where $\mu(n)$ is the M\"obius function.

Using this formula and the Poincar\'e-Birkhoff-Witt theorem
$\Ass^+=\Com^+\circ\Lie^+$, we may calculate $\ch(\Lie^+)$.
\begin{proposition} \label{Lie}
$\displaystyle\ch(\Lie^+) = - \sum_{n=1}^\infty \frac{\mu(n)}{n}
\log(1-p_n)$
\end{proposition}
\begin{pf}
We know from the Poincar\'e-Birkhoff-Witt theorem that
$\ch(\Ass^+)=\Exp(\ch(\Lie^+))$; it follows that
$$
\ch(\Lie^+) = \Log(1+\ch(\Ass^+))
= \sum_{n=1}^\infty \frac{\mu(n)}{n} \log(1+p_n\circ\ch(\Ass^+)) .
$$
Since
$$
1+ p_n \circ \ch(\Ass^+) = 1 + \frac{p_n}{1-p_n} = \frac{1}{1-p_n} ,
$$
the result follows.
\end{pf}

It follows from this formula that
$$
\ch_n(\Lie^+(n)) = \frac{1}{n} \sum_{d|n} \mu(d) p_d^{n/d} .
$$
This is the characteristic of the induced representation
$\Ind_{C_n}^{\SS_n}\chi$, where $\chi$ is a primitive character of the
cyclic group.

We now turn to calculating the characteristic of the braid operad $\Braid$.
First, we need a lemma.
\begin{lemma}
$\ch_t(\Lambda\v)=-t\ch_t(\v)(-t^{-1}p_1,-t^{-2}p_2,-t^{-3}p_3,\dots)$
\end{lemma}
\begin{pf}
Tensoring with $\sgn_n$ has the effect of replacing $p_n$ by
$(-1)^{n-1}p_n$. Applying $\Sigma^{-n}$ to $\v(n)$ then has the effect of
replacing $p_n$ by $(-t)^{-n}p_n$.
\end{pf}

\begin{proposition} \label{configuration}
For each $n\ge1$, let
$$
P_n(t) = \frac{1}{n} \sum_{d|n} \frac{\mu(n/d)}{t^d} .
$$
Then
$$
\ch_t(\Braid) = \prod_{n=1}^\infty \Bigl( 1 + t^np_n \Bigr)^{P_n(t)} - 1 .
$$
\end{proposition}
\begin{pf}
The $\SS$-module $\Lambda^{-1}\Lie^+$ has Poincar\'e characteristic
$$
\ch_t(\Lambda^{-1}\Lie^+)
= \frac{1}{t} \sum_{n=1}^\infty \frac{\mu(n)}{n} \log(1+t^np_n) .
$$
It follows that
\begin{eqnarray*}
\ch_t(\b) &=& \Exp \circ \ch_t(\Lambda^{-1}\Lie^+) - 1
= \exp\Bigl( \sum_{k=1}^\infty \frac{p_k}{k} \Bigr) \circ
\Bigl( t^{-1} \sum_{n=1}^\infty \frac{\mu(n)}{n} \log(1+t^np_n) \Bigr) - 1 \\
&=& \exp\Bigl( \sum_{k=1}^\infty
\sum_{n=1}^\infty t^{-k} \frac{\mu(n)}{kn} \log(1+t^{kn}p_{kn}) \Bigr) - 1
\\
&=& \prod_{n=1}^\infty ( 1 + t^np_n )^{P_n(t)} - 1 .
\qed
\end{eqnarray*}
\def\qed{}\end{pf}

In particular, setting $p_1=x$ and $p_n=0$, $n>1$, we see that the
Poincar\'e polynomial of the space $\C^n_0$ is
\begin{eqnarray}
\label{poincare-config}
\sum_{i=0}^{n-1} (-t)^n \dim H^i(\C^n_0)
&=& \text{coefficient of $x^n/n!$ in $(1+tx)^{t^{-1}} - 1$} \\
&=& \binom{t^{-1}}{n} t^n = \prod_{i=1}^{n-1} (1-it) .
\end{eqnarray}

\subsection{The characteristic of a cyclic $\SS$-module}
If $\v$ is a stable cyclic $\SS$-module, we define $\Ch_t(\v)$ in a similar
way to $\ch_t(\v)$:
$$
\Ch_t(\v) = \sum_{n=3}^\infty \sum_{i=0}^\infty (-t)^i \,
\ch_n(\v_i\(n\)) \in \Lambda\(t\) .
$$
If $\ch_t(\v)$ denotes the Poincar\'e characteristic of the $\SS$-module
underlying $\v$, then we have the formula
$$
\ch_t(\v) = \frac{\p\Ch_t(\v)}{\p p_1} .
$$

We now calculate $\Ch_t(\m)$, where $\m$ is the cyclic $\SS$-module
$$
\m\(n\) = \begin{cases} H_\bull(\CM_{0,n}) , & n\ge3 , \\ 0 , & n<3 .
\end{cases}$$
Note that the Poincar\'e polynomial of $\CM_{0,n}$ is much easier to
calculate than the Poincar\'e characteristic: it is obtained by dividing
the Poincar\'e polynomial \eqref{poincare-config} of $\C^{n-1}_0$ by $1-t$:
\begin{equation} \label{poincare-moduli}
\sum_{i=0}^{n-2} (-t)^n \dim H^i(\CM_{0,n}) = \prod_{i=2}^{n-2} (1-it).
\end{equation}
\begin{theorem} \label{moduli:characteristic}
$$
\Ch_t(\m) = \frac{1}{1-t^2} \Bigl( (1+tp_1)
\prod_{n=1}^\infty ( 1 + t^np_n )^{P_n(t)}
- 1 - (1+t)h_1 - (h_2 + te_2) \Bigr)
$$
\end{theorem}
\begin{pf}
{}From the $\SS_n$-equivariant homotopy equivalence
$\C^n_0\simeq\CM_{0,n+1}\times S^1$, which holds for $n\ge2$, we see that
\begin{equation} \label{differential}
\frac{\p\Ch_t(\m)}{\p p_1}
= \ch_t(\m) = \frac{t\bigl(\ch_t(\Braid)-p_1\bigr)}{t-1} .
\end{equation}
The Serre spectral sequence for the $\SS_n$-equivariant fibration
$$\begin{CD}
\C\setminus\{1,\dots,n\} @>>> \CM_{0,n+1} @>>> \CM_{0,n}
\end{CD}$$
collapses at $E^2$, so $H_\bull(\CM_{0,n+1}) \cong
H_\bull(\C\setminus\{1,\dots,n\}) \o H_\bull(\CM_{0,n})$.
Furthermore, this isomorphism is $\SS_n$-equivariant, where $\SS_n$
acts on $H_\bull(\C\setminus\{1,\dots,n\})$ by the monodromy of the
Gauss-Manin connection.

Now, $H_0(\C\setminus\{1,\dots,n\})$ is the trivial $\SS_n$-module, while
$H_1(\C\setminus\{1,\dots,n\})$ is the irreducible representation
$V_{n-1,1}$, which is the kernel of the natural map $\begin{CD}
\C^n@>>>\C\end{CD}$ obtained
by sending $(x_1,\dots,x_n)$ to $x_1+\dots+x_n$. If $\sigma\in\SS_n$ is a
transitive permutation, $\Tr(\sigma|V_{n-1,1})=-1$; this shows that
\begin{equation} \label{initial}
\Ch_t(\m)\big|_{p_1=0} = \frac{\ch_t(\Braid)}{1+t} \Big|_{p_1=0} .
\end{equation}
The theorem now follows on solving the differential equation
\eqref{differential} with initial condition \eqref{initial}.
\end{pf}

The first few terms of $\Ch_t(\m)$ are as follows:
$$\begin{tabular}{|c|l|}
$n$ & $\Ch_t(\CM_{0,n})$ \\ \hline
$3$ & $s_3$ \\
$4$ & $s_4-ts_{2^2}$ \\
$5$ & $s_5-ts_{32}+t^2(s_{31^2})$ \\
$6$ & $s_6-ts_{42}+t^2(s_{41^2}+s_{321})-t^3(s_{41^2}+s_{3^2}+s_{2^21^2})$ \\
\end{tabular}$$
The pattern emerging here, that $H_1(\CM_{0,n})\cong V_{n-2,2}$, is easily
verified in general using our formula for $\Ch_t(\m)$. We have seen that
there is a natural identification between $\sgn_n\o H^1(\CM_{0,n})$ and the
space of relations \eqref{hypercommutative} among $n-1$ letters in the
hypercommutative operad: thus, we see that this space of relations is the
irreducible $\SS_n$-module $V_{2^21^{n-4}}$.

Applying l'H\^opital's rule to Theorem \ref{moduli:characteristic}, we see
that the Euler-Frobenius characteristic of $\m$ is given by the formula
$$
\Ch(\m) = \lim_{\begin{CD}t@>>>1\end{CD}} \Ch_t(\m)
= \frac12 (1+p_1)^2 \sum_{n=1}^\infty \frac{\varphi(n)}{n}
\log(1+p_n) - \frac{1}{4} (2p_1 + 3p_1^2 + p_2) .
$$

Finally, it follow easily from the formula for $\Ch_t(\m)$ that
$$
\Ch_t(\Grav) = - \frac{t^3}{1-t^2}
\Bigl( (1-p_1) \prod_{n=1}^\infty (1-p_n)^{P_n(t)}
- 1 + (1+t^{-1})h_1 - (t^{-1}h_2 + t^{-2}e_2) \Bigr)
$$

\subsection{The Poincar\'e characteristic of $\TT_+\v$}
In \cite{modular}, a formula for the Poincar\'e characteristic of $\TT_+\v$
in terms of the Poincar\'e characteristic of $\v$ is derived. If
$F=e_2-\Ch_t(\v)$, define the Legendre transform $G=\CL F$ of $F$ in the
sense of symmetric functions by the formula
\begin{equation} \label{Legendre}
F\circ \frac{\p(\CL G)}{\p p_1} + G = p_1 \frac{\p G}{\p p_1} .
\end{equation}
Then $G=h_2+\Ch_t(\TT_+\v)$.

Note that \eqref{Legendre} implies that $(\p F/\p p_1)\circ(\p G/p p_1)=p_1$,
from which it is straightforward to calculate $\p G/\p p_1$. Substituting
$\p G/\p p_1$ into both sides of \eqref{Legendre}, we obtain an explicit
formula for $G$.

As an application of \eqref{Legendre}, we now calculate the Poincar\'e
characteristics of the varieties $\Mbar_{0,n}$. We use a slight extension
of \eqref{Legendre}, in which $\v=\v_0\oplus\v_1$ has an internal
$\Z/2$-grading, and $\Ch_t(\v)=\Ch_t(\v_0)-\Ch_t(\v_1)$. Let $\v$ and $\w$
be the stable cyclic $\SS$-modules
\begin{eqnarray*}
\v_i\(n\) &=& \begin{cases} 0 , & i=0 , \\ \mbar\(n\)\o\sgn_n , & i=1 ;
\end{cases} \\
\w_{ij}\(n\) &=& \begin{cases} \sgn_n\o H^p(\CM_{0,n}) , &
\text{$i=2(n-p-3)$ and $j\equiv p+1\pmod{2}$,} \\
0 , & \text{otherwise.}
\end{cases}
\end{eqnarray*}
Thus,
\begin{eqnarray*}
\Ch_t(\v) & =& - \Ch_t(\mbar)\big|_{p_n\to(-1)^{n-1}p_n} , \\
\Ch_t(\w) &=& \frac{- t^{-6}\Ch_t(\m)\Big|_{t\to-t^{-2}}
{p_n\to(-1)^{n-1}t^{2n}p_n}} .
\end{eqnarray*}
Using \eqref{det}, we may rewrite \eqref{bar:even} in the form
$$\begin{CD}
0 @>>> H^p(\CM_{0,n})\o\sgn_n @>>>
\bigoplus_{T\in\CT_p\(n\)} \v\(T\)_{2(n-p-3)} @>>>
\bigoplus_{T\in\CT_{p-1}\(n\)} \v\(T\)_{2(n-p-3)} @>>> \dots
\end{CD}$$
which shows that $\Ch_t(\w)=\Ch_t(\TT_+\v)$, and hence that
$G=h_2+\Ch_t(\w)$ is the Legendre transform of $F=e_2-\Ch_t(\v)$. In this
way, we have proved the following proposition.
\begin{theorem} \label{F-G}
The symmetric function
$$
F = e_2 + \Ch_t(\mbar)\big|_{p_n\to(-1)^{n-1}p_n}
$$
is the Legendre transform of the symmetric function
$$
G = h_2 - t^{-6} \Ch_t(\m)\Big|_{\frac{t\to t^2}{ p_n\to(-1)^{n-1}t^{2n}p_n}} .
$$
\end{theorem}

Here are some sample calculations of $\Ch_t(\Mbar_{0,n})$ for small $n$:
$$\begin{tabular}{|c|l|}
$n$ & $\Ch_t(\Mbar_{0,n})$ \\ \hline
$3$ & $s_3$ \\
$4$ & $(1+t^2)s_4$ \\
$5$ & $(1+t^4)s_5+t^2(s_5+s_{41})$ \\
$6$ & $(1+t^6)s_6+(t^2+t^4)(2s_6+s_{51}+s_{42})$ \\
\end{tabular}$$

These formulas simplify if we are only interested in the dimensions of
the vector spaces $H_i(\Mbar_{0,n})$. We have the formula
\begin{eqnarray*}
g(x,t) = G'\Big|_{\frac{p_1\to x}{ p_n\to0,n>1}}
&=& x - \sum_{n=2}^\infty \frac{x^n}{n!} \sum_{i=0}^{n-2} (-1)^i
t^{2(n-i-2)} \dim H_i(\CM_{0,n+1}) \\
&=& x - \frac{(1+x)^{t^2}-(1+t^2x)}{t^2(t^2-1)} .
\end{eqnarray*}
It is a corollary of Theorem \ref{F-G} that
$$
f(x,t) = F'\Big|_{\frac{p_1\to x}{ p_n\to0,n>1}}
= x + \sum_{n=2}^\infty \frac{x^n}{n!}
\sum_{i=0}^{n-2} t^{2i} \dim  H_{2i}(\Mbar_{0,n+1}) ,
$$
is the inverse of $g$, in the sense that $f(g(x,t),t)=x$.
This is a reformulation of Fulton and MacPherson's calculation of the
Poincar\'e polynomial of $\Mbar_{0,n}$. Note that their proof also makes
use of mixed Hodge theory, in the form of the ``fake Poincar\'e
polynomial.'' Our result Theorem \ref{F-G} is an equivariant version of their
calculation.

\makeatletter\renewcommand{\@biblabel}[1]{\hfill[#1]}\makeatother


\begin{thebibliography}{99}

\bibitem{AC} E. Arbarello and M. Cornalba, \emph{The Picard groups of the
moduli spaces of curves,} Topology \textbf{26} (1987), 153--171.

\bibitem{Arnold}
V.I. Arnold, \emph{The cohomology ring of the colored braid group,}
Mat. Zametki \textbf{5} (1969), 227--231.

\bibitem{BG} A. Beilinson, V. Ginzburg, \emph{Infinitesimal structure of
moduli spaces of $G$-bundles,} Internat. Math. Res. Notices (appendix
to Duke Math. J.) \textbf{66} (1992), 63--74.

\bibitem{BoardmanVogt}
J.M. Boardman and R.M. Vogt, ``Homotopy invariant algebraic structures on
topological spaces,'' Lecture Notes in Math. \textbf{347}, 1973.

\bibitem{BZ} J.L. Brylinski and S. Zucker, \emph{An overview of recent
advances in Hodge theory,} in ``Several complex variables, VI,'' 39--142,
Encyclopaedia Math. Sci. \textbf{69}, Springer Verlag, Berlin, 1990.

\bibitem{Cohen:thesis}
F.R. Cohen, \emph{The homology of {$\mathcal{C}_{n+1}$}-spaces,
{$n\ge0$}}, in ``The homology of iterated loop spaces,'' Lecture
Notes in Math. \textbf{533}, 1976, 207--351.

\bibitem{Deligne} P. Deligne, \emph{Th\'eorie de Hodge II,} Publ. Math.
IHES \textbf{40} (1971), 5--58.

\bibitem{DVV} R. Dijkgraaf, E. and H. Verlinde, \emph{Topological strings
in $d<1$,} Nucl. Phys. \textbf{B352} (1991), 59--80.

\bibitem{FM}
W. Fulton and R. MacPherson, \emph{A compactification of configuration
spaces}, Ann. Math., \textbf{139} (1994), 183--225.

\bibitem{bv}
E. Getzler, \emph{Batalin-Vilkovisky algebras and two-dimensional
topological field theories,} Commun. Math. Phys. \textbf{159} (1994),
265--285.

\bibitem{weil} E. Getzler, \emph{Equivariant cohomology and topological
gravity,} Commun. Math. Phys. \textbf{163} (1994), 473--490.

\bibitem{n-algebras}
E. Getzler and J.D.S. Jones, \emph{Operads, homotopy algebra, and
iterated integrals for double loop spaces,} (\texttt{hep-th/9403055}).

\bibitem{cyclic}
E. Getzler and M. Kapranov, \emph{Cyclic operads and cyclic homology,} to
appear in ``Geometry, Topology, and Physics for Raoul,'' ed. B. Mazur,
International Press, Cambridge, MA, 1994.

\bibitem{modular}
E. Getzler and M. Kapranov, \emph{Modular operads,} MPIM-Bonn preprint
94/78, (\texttt{dg-ga/9408003}).

\bibitem{GK}
V.A. Ginzburg and M.M. Kapranov, \emph{Koszul duality for operads,}
to appear, Duke. Math. J.

\bibitem{Gromov} M. Gromov, \emph{Pseudo holomorphic curves in
symplectic manifolds,} Invent. Math. \textbf{82} (1985), 307--347.

\bibitem{Joyal}
A. Joyal, \emph{Foncteurs analytiques et esp\`eces de structures,} Lecture
Notes in Math. \textbf{1234} (1986) 126--159.

\bibitem{Kapranov} M.M. Kapranov, \emph{Permuto-associahedron, MacLane's
coherence theorem and asymptotic zones for the KZ equation,} J. Pure Appl.
Algebra, \textbf{85} (1993), 119--142

\bibitem{Keel}
S. Keel, \emph{Intersection theory of moduli spaces of stable $n$-pointed
curves of genus zero,} Trans. Amer. Math. Soc. \textbf{330} (1992),
545--574.

\bibitem{Klyachko} A.A. Klyachko, \emph{Lie elements in the tensor algebra,}
Siberian Math. J., \textbf{15} (1974), 914--920.

\bibitem{Knudsen} F.F. Knudsen, \emph{The projectivity of the moduli space
of stable curves II. The stacks $\Mbar_{g,n}$,} Math. Scand. \textbf{52}
(1983), 161--189.

\bibitem{Kontsevich}
M. Kontsevich, \emph{Formal (non)-commutative symplectic geometry,\/}
in ``The Gelfand mathematics seminars, 1990--1992,'' eds. L. Corwin,
I. Gelfand, J. Lepowsky, Birkh\"auser, Boston, 1993.

\bibitem{KM}
M. Kontsevich, Yu. Manin, \emph{Gromov-Witten classes, quantum cohomology,
and enumerative geometry,} Commun. Math. Phys. \textbf{164} (1994).

\bibitem{Macdonald}
I.G. Macdonald, ``Symmetric Functions and Hall Polynomials,'' Clarendon
Press, Oxford, 1979.

\bibitem{MS}
D. McDuff and D. Salamon, ``$J$-holomorphic curves and quantum
cohomology,'' Amer. Math. Soc., Providence, 1994.

\bibitem{Maclane}
S. Maclane, ``Categories for the working mathematician,'' Graduate
Texts in Math. \textbf{5}, 1971.

\bibitem{Quillen}
D. Quillen, \emph{Rational homotopy theory,} Ann. Math. \textbf{90}
(1969), 205--295.

\bibitem{RT} Y. Ruan, G. Tian, \emph{A mathematical theory of quantum
cohomology,} preprint, 1994.

\end{thebibliography}
\end{document}